\documentclass[reqno, 11pt]{article}
\usepackage{amssymb, amsmath, amsthm, amsfonts, amscd, epsfig,bm}
\usepackage{color}
\usepackage[english]{babel}
\usepackage{bbm}
\usepackage{geometry,graphicx,pgfplots,caption,subcaption, diagbox}
\usepackage[export]{adjustbox}
\usepackage[titletoc,toc,title]{appendix}
\usepackage[normalem]{ulem} 

\usepackage[colorlinks=true,
            linkcolor=blue,
            citecolor=blue,
            urlcolor=blue]{hyperref}

\newtheorem{thm}{Theorem}[section]

\newtheorem{lemma}[thm]{Lemma}

{\bf}{\it}

\def\N{\mathbb{N}}
\def\Z{\mathbb{Z}}

\def\C{\mathbb{C}}

\newcommand{\Om}{\Omega}

\newcommand{\RR}{\mathbb{R}}

\newcommand{\eps}{\varepsilon}

\newcommand{\p}{\partial}

\newcommand{\eqnref}[1]{(\ref {#1})}

\newcommand{\beq}{\begin{equation}}
\newcommand{\eeq}{\end{equation}}
\newcommand{\RN}[1]{%
\textup{\uppercase\expandafter{\romannumeral#1}}%
}

\numberwithin{equation}{section}
\numberwithin{figure}{section}

\begin{document}







\title{
Convergence analysis of resonant frequency approximation in periodic Helmholtz problems
}

\date{}

\author{
Jiho Hong\thanks{Department of Mathematics, The Chinese University of Hong Kong, Shatin, N.T., Hong Kong, P.R. China (jihohong@cuhk.edu.hk).}\and
Woojoo Lee\thanks{Department of Mathematical Sciences, Korea Advanced Institute of Science and Technology, 291 Daehak-ro, Yuseong-gu, Daejeon 34141, Republic of Korea (woojoo.lee@kaist.ac.kr, mklim@kaist.ac.kr).}\and
Mikyoung Lim\footnotemark[3]
}

\maketitle

\begin{abstract}
We investigate resonant frequencies of the Helmholtz equation in periodic structures under quasi-periodic conditions, using approximations of the associated quasi-periodic Green function, which is represented by a lattice sum. By formulating the periodic Helmholtz problem using layer potential operators whose kernel is given by the quasi-periodic Green function, resonant frequencies correspond to characteristic values of an associated operator-valued function. 
Based on two approximation methods for the quasi-periodic Green function, the Ewald method and the shifted Green function method, we propose the corresponding approximation schemes for resonant frequencies. Furthermore, we establish convergence rates for the resulting approximations using Gohberg--Sigal theory, where the Ewald method yields exponential convergence, while the shifted Green function method yields algebraic convergence. Numerical computations are presented to validate and compare the proposed approximation schemes.
\end{abstract}


\tableofcontents


\section{Introduction}
Periodic structures play an important role in various wave-scattering applications, including metascreens composed of subwavelength plasmonic particles \cite{Ammari:2016:MNF} and Helmholtz resonators \cite{Ammari:2015:MTS}. Related wave phenomena in resonant media also have applications in ultrasonic imaging and focusing \cite{Errico:2015:UUL, Lanoy:2015:SFB}.  
In particular, bubbly media exhibit distinctive wave-propagation properties even when the bubbles occupy only a small volume fraction, leading to changes in the effective wave speed compared with that of the surrounding fluid \cite{Ammari:2017:EMT, Commander:1989:LPW, Hwang:2000:LFR}.
In the subwavelength regime, where the wavelength is much larger than the bubble size, bubbles can exhibit Minnaert resonance \cite{Minnaert:1933:MAB}; see \cite{Ammari:2018:MRA, Calvo:2012:LFS, Leroy:2015:SAW} for related studies. Such resonant phenomena have been exploited to achieve effects such as the superabsorption of acoustic waves \cite{Leroy:2015:SAW,Leroy:2009:TUS} and have been investigated in connection with the formation of band gaps in phononic crystals \cite{Leroy:2009:DCB}. 
Further approaches have been developed for periodic structures in optics \cite{Darweesh:2018:RRW,Maradudin:2016:RWA,Savoia:2013:SSR,Tavakoli:2019:RWA}, as well as in computational acoustic and electromagnetic problems \cite{Cho:2019:SAN,Lai:2015:FRS,Liu:2016:ENS}.

The analysis of wave scattering by periodic structures may involve quasi-periodic conditions for the incident fields and scattered waves. Resonant frequencies in periodic acoustic structures have been investigated both analytically and numerically. For instance, \cite{Ammari:2017:MNF} provides a mathematical and numerical framework for bubbly meta-screens consisting of periodic subwavelength bubbles mounted on a reflective surface. Related analytical approaches to acoustic resonances, including Minnaert resonances of bubbles of arbitrary shapes, are developed in \cite{Ammari:2018:MRA}. Numerical results in these works also illustrate how resonant frequencies depend on the geometric configuration of the bubbles \cite{Ammari:2015:MTS, Ammari:2016:MNF,Ammari:2017:MNF}. In the context of electromagnetic scattering, plasmonic resonances have also been extensively studied and characterized through the spectrum of Neumann--Poincaré integral operators \cite{Petit:1990:ETG}. 

Boundary integral formulations are widely used to analyze wave scattering by periodic structures, in which the quasi-periodic Green function plays a central role \cite{Ammari:2009:LPT}. However, its classical lattice-sum representation converges slowly, making its direct numerical evaluation computationally demanding. Various techniques have therefore been developed to accelerate the convergence based on its alternative representations, including the Laplace-type integral method, the Ewald method \cite{Ewald:1921:BOE}, and the shifted Green function method \cite{Bruno:2014:RCT}. 
the Ewald method has often been employed to obtain rapidly convergent representations \cite{Capolino:2007:EC3}, while the shifted Green function method provides another approach to efficient evaluation \cite{Bruno:2014:RCT,Bruno:2017:RCQ,Bruno:2017:TDQ}. In particular, the latter can yield rapidly convergent representations even in the presence of Rayleigh–Wood anomalies, where the classical quasi-periodic Green function may cease to exist \cite{Maystre:2012:TWA}. For alternative approaches regarding domain decomposition, we refer to \cite{Perez-Arancibia:2019:DDQ, Nicholls:2020:SPI}.

In this paper, we investigate approximation methods for resonant frequencies in periodic structures. Focusing on the case of a single inclusion per periodic cell, we formulate the resonant frequencies as characteristic values of an operator-valued function and develop approximation schemes. While previous studies have compared various methods for approximating the quasi-periodic Green function \cite{Bruno:2020:EQG}, we consider two of these methods: the Ewald method and the shifted Green function method. Based on these approximations, we employ Gohberg--Sigal theory (see Appendix \ref{appen:GS}) to establish convergence of the approximated characteristic values and show that their convergence rates are governed by the corresponding approximation of the quasi-periodic Green function.

Using layer potential operators, the Helmholtz problem in a periodic domain can be formulated as a boundary integral equation involving an operator-valued function, whose characteristic values correspond to the resonant frequencies \cite{Ammari:2009:LPT}. The key idea of our analysis is to transfer the convergence of the quasi-periodic Green function approximations to that of the associated characteristic values. To this end, we employ Gohberg--Sigal theory, which provides a framework for analyzing the stability of characteristic values under perturbations of operator-valued functions. By approximating the boundary integral operator through truncated representations of the quasi-periodic Green function, we derive convergence rates for the resulting resonant frequency approximations. In particular, we establish exponential convergence for the Ewald method and algebraic convergence for the shifted Green function method. Numerical experiments further illustrate the convergence behavior of the proposed schemes, complementing the theoretical convergence analysis.

The rest of the paper is structured as follows. In Section \ref{sec:prob}, we introduce the periodic Helmholtz problem and formulate the problem using layer potential operators, which involve the quasi-periodic Green function. We review two approximation methods for the quasi-periodic Green function and derive their convergence rates in Section \ref{sec:Green}.
Based on these results, in Section \ref{sec:resonant}, we propose a computational scheme for determining resonant frequencies in periodic domains and present its convergence behavior. The proof of the corresponding theorem is provided in Section \ref{subsec:resonant:proof}. We validate our findings through numerical computation in Section \ref{sec:num} and conclude the paper in Section \ref{sec:conclusion}.  

\section{Problem formulation} \label{sec:prob}
\subsection{Periodic Helmholtz problems} \label{subsec:prob:periodic} 
We first introduce a periodic structure and state the related scattering problems.
Let $d=2$ or $3$ and $\Lambda$ be a $d_\Lambda$-dimensional periodic lattice in $\RR^d$, where $1\le d_{\Lambda} < d$. For any vector $\mathbf{x} \in \RR^d$, we denote the projection of $\mathbf{x}$ onto $\text{span}(\Lambda)$ by $\mathbf{x}^\parallel$ and the projection of $\mathbf{x}$ onto $\text{span}(\Lambda)^\perp$ by $\mathbf{x}^\perp$. Let $\Om_p$ ($0\le p\le P$) be the $\Lambda$-periodic disjoint domains in $\mathbb{R}^d$ such that $\p\Om_p$ are smooth, and let the refractive index be piecewise constant, with $n(\mathbf{x})=\sum_{p=0}^{P}n_p\chi_{\Om_p}(\mathbf{x})$.
We fix vectors $\mathbf{\alpha}$ and $\mathbf{\beta} \in \mathbb{R}^d$ satisfying
$$|\mathbf{\alpha}|^2+|\mathbf{\beta}|^2=n^2 \omega^2,\quad \mathbf{\alpha}\parallel \Lambda\ \ \mbox{and}\ \ \mathbf{\beta}\perp\Lambda.$$
In the presence of an incident field $u^{\rm inc}(\mathbf{x})=e^{{\rm i}\mathbf{\alpha}\cdot\mathbf{x}-{\rm i}\mathbf{\beta}\cdot\mathbf{x}}$, we assume the total field $u$ satisfying
\begin{align}
&\nabla\cdot\left(\frac{1}{n(\mathbf{x})^2}\nabla u\right) +  \omega^2 u =0 \quad\mbox{in }\Om_p
\end{align}
in the sound-hard case, and
\begin{align}
&\Delta u +  \omega^2{n(\mathbf{x})^2} u =0 \quad\mbox{in }\Om_p
\end{align}
in the sound-soft case, where the scattered field $u^{\rm s}:=u-u^{\rm inc}$ is subject to the radiation condition at infinity described in \cite[Eq. (21-22)]{Bruno:2020:EQG}.
On the interfaces $\Gamma_{pq}$ between adjacent domains $\Omega_p$ and $\Omega_q$, there hold the transmission conditions
\begin{align*}
u_p|_{\Gamma_{pq}} =  u_q|_{\Gamma_{pq}}\quad\mbox{and}\quad \p_\nu u_p|_{\Gamma_{pq}} = C_{pq}^2\p_\nu u_q|_{\Gamma_{pq}},
\end{align*}
where we set $u_p=u|_{\Om_p}$ for $p=0,1,\dots,P$, $\nu$ to be the unit normal vector to $\Gamma_{pq}$, and $C_{pq}=n_p/n_q$ in the sound-hard case and $C_{pq}=1$ in the sound-soft case.
The incident field is $\mathbf{\alpha}$-quasi-periodic in the sense that $e^{-{\rm i}\mathbf{\alpha} \cdot\mathbf{x}}u^{\rm inc}$ is $\Lambda$-periodic.
It is well-known \cite{Petit:1990:ETG} that $u^{\rm s}$ is also $\mathbf{\alpha}$-quasi-periodic.
For the comparison of various computation schemes, we refer the reader to \cite{Bruno:2020:EQG}. 

\begin{figure}[t]
	\centering
	\includegraphics[width=0.5\textwidth]{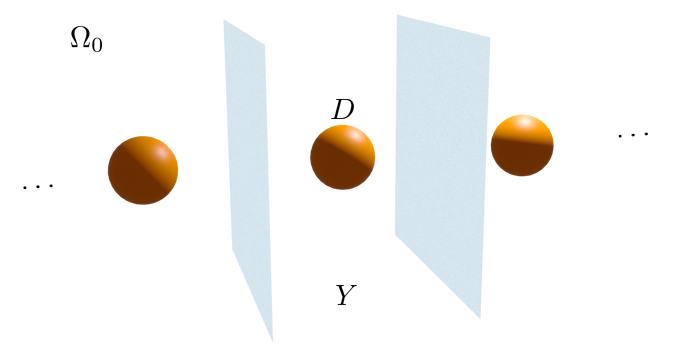} 
	\caption{\label{fig:periodic} An example of a periodic structure $\Omega_1 = \cup_{\lambda \in \Lambda} (\lambda + D)$ with $d=3$ and $d_\Lambda=1$, illustrating the primitive cell $Y$ and the bounded domain $D \subset Y$.} 
\end{figure}

We consider resonances in the periodic structures with a single inclusion per cell. 
For simplicity, we restrict our attention to the case when $P=1$ with $\Omega_1=\cup_{\lambda\in\Lambda}(\lambda+D)$ for a bounded domain $D$ in the primitive cell of $\Lambda$ and $\Omega_0=\mathbb{R}^d\backslash \overline{\Omega_1}$.
In other words, we seek for the eigenvalue $\omega^2$ of the problem
\begin{equation}\label{eq:eigproblem:simple}
\left\{\begin{aligned}
\Delta u + n_0^2\omega^2 u =0&\quad\mbox{in }Y\backslash\overline{D},\\
\Delta u + n_1^2\omega^2 u =0&\quad\mbox{in }{D},\\
u|_+ =u|_-&\quad\mbox{on }\partial{D},\\
\partial_\nu u|_+ =C_{01}^2\partial_\nu u|_- &\quad\mbox{on }\partial{D},\\
e^{-{\rm i}\mathbf{\alpha}\cdot\mathbf{x}} u\mbox{ is }\Lambda\mbox{-periodic},
\end{aligned}\right.
\end{equation}
where $Y$ is the primitive cell of the lattice $\Lambda$ (see \ref{fig:periodic}).

\subsection{Layer potential formulation} \label{subsec:prob:lp} 
Let $\{\mathbf{R}_m\}_{m \in \mathbb{Z}^{d_\Lambda}}$ be an enumeration of $\Lambda$, so that $\Lambda=\{\mathbf{R}_m\,:\,m\in\mathbb{Z}^{d_\Lambda}\}$. The quasi-periodic Green function associated with $\Lambda$ is given by
\begin{equation} \label{eq:Gq}
G_\kappa^{q}(\mathbf{x})=\sum_{m\in\mathbb{Z}^{d_\Lambda}}e^{{\rm i}{\bf\alpha}\cdot\mathbf{R}_m}G_\kappa(\mathbf{x}-\mathbf{R}_m),
\end{equation}
where $G_\kappa$ is the free-space Green function for the Helmholtz equation, defined by
$$G_\kappa(\mathbf{x})=\begin{cases}
\frac{\rm i}{4}H_0^{(1)}(\kappa|\mathbf{x}|)&\mbox{for }d=2,\\
\frac{1}{4\pi}\frac{e^{{\rm i}\kappa|\mathbf{x}|}}{|\mathbf{x}|}&\mbox{for }d=3.
\end{cases}$$
The corresponding single- and double-layer potential operators are given by
\begin{align} 
\label{def:S}&\mathcal{S}_{\Om}^{q}(\kappa)[\varphi](\mathbf{x})=\int_{\p\Om} G_\kappa^q(\mathbf{x}-\mathbf{y})\varphi(\mathbf{y})\,{\rm d}\sigma(\mathbf{y}),\quad\varphi\in L^2(\p\Om),\\
&\mathcal{D}_{\Om}^{q}(\kappa)[\varphi](\mathbf{x})=\int_{\p\Om} \frac{\p}{\p\nu(\mathbf{y})}G_\kappa^q(\mathbf{x}-\mathbf{y})\varphi(\mathbf{y})\,{\rm d}\sigma(\mathbf{y}),\quad\varphi\in L^2(\p\Om),
\end{align}
where $\nu$ is the unit normal vector pointing outward of $\Om$.
We also consider the operators
\begin{align}
&\mathcal{K}_{\Om}^{q}(\kappa)[\varphi](\mathbf{x})=\int_{\p\Om} \frac{\p}{\p\nu(\mathbf{y})}G_\kappa^q(\mathbf{x}-\mathbf{y})\varphi(\mathbf{y})\,{\rm d}\sigma(\mathbf{y}),\quad\varphi\in L^2(\p\Om),\\
\label{def:Kstar}&\mathcal{K}_{\Om}^{q,*}(\kappa)[\varphi](\mathbf{x})=\int_{\p\Om} \frac{\p}{\p\nu(\mathbf{x})}G_\kappa^q(\mathbf{x}-\mathbf{y})\varphi(\mathbf{y})\,{\rm d}\sigma(\mathbf{y}),\quad\varphi\in L^2(\p\Om).
\end{align}

We recall some properties of layer potential operators from \cite[Chapter 8]{Ammari:2009:LPT}. 
\begin{lemma}[{\cite[Theorem 8.12]{Ammari:2009:LPT}}]
Let $u$ be the eigenfunction of \eqref{eq:eigproblem:simple}. If $n_0\omega$ and $n_1\omega$ are not Dirichlet eigenvalues in $Y\backslash\overline{D}$ and $D$, then there is only one pair $(\varphi,\psi)\in L^2(\p D)^2$ such that the following representation holds:
\begin{equation}\label{eq:repr:eigenfunction}
u=\begin{cases}
\mathcal{S}_{D}^{q}(n_1\omega)\varphi&\mbox{in }D,\\
\mathcal{S}_{D}^{q}(n_0\omega)\psi&\mbox{in }Y\backslash\overline{D}.
\end{cases}
\end{equation}
\end{lemma}
Under the representation \eqref{eq:repr:eigenfunction}, the problem \eqref{eq:eigproblem:simple} is equivalent to the following system of boundary integral equations:
\begin{align} \label{def:A}
\mathcal{A}(\omega)\begin{pmatrix}
\varphi\\ \psi    
\end{pmatrix}=0,\quad\mbox{where}\quad \mathcal{A}(\omega)=\begin{pmatrix}
\mathcal{S}_{D}^{q}(n_1\omega) & -\mathcal{S}_{D}^{q}(n_0\omega)
\\
C_{01}^2\left(-\frac{1}{2}{\rm id}+\mathcal{K}_{D}^{q,*}(n_1\omega)\right)&-\frac{1}{2}{\rm id}-\mathcal{K}_{D}^{q,*}(n_0\omega)
\end{pmatrix},
\end{align}
where $C_{01} = n_0/n_1$, and ${\rm id}$ denotes the identity operator. Then the resonant frequency $\omega$ corresponds to the characteristic value of the operator-valued function $\mathcal{A}(\omega)$.
We further have the following characterization of the resonant frequencies, as a consequence of Theorem \ref{thm:GS:Steinberg}:
\begin{lemma}[{\cite[Lemma 8.1]{Ammari:2009:LPT}}]
\label{lemma:poles:characterization}
The operator-valued function $\mathcal{A}(\omega)$ is Fredholm analytic with index $0$ in $\mathbb{C}\backslash{\rm i}\mathbb{R}^-$. Moreover, $\omega\mapsto (\mathcal{A}(\omega))^{-1}$ is meromorphic and its poles are on the real axis.
\end{lemma}

\section{Approximation of quasi-periodic Green function} \label{sec:Green}

We review two approximation methods for the quasi-periodic Green function \eqref{eq:Gq}, namely, the Ewald method and the shifted Green function method. We introduce the methods and obtain their convergence rates. 

In this work, we fix $d=2$ or $d=3$ with $1 \leq d_\Lambda <d$. We denote the primitive vectors of $\Lambda$ by $\mathbf{v}_j$ for $j=1,\ldots,d_\Lambda$. By applying a suitable change of basis, we assume that $\text{span}(\{\mathbf{v}_1, \ldots, \mathbf{v}_{d_\Lambda} \}) = \text{span}(\{\hat{\mathbf{e}}_1, \ldots, \hat{\mathbf{e}}_{d_\Lambda}\})$. Let $A$ be the area of the primitive cell, defined by 
$$A = \begin{cases}
|\mathbf{v}_1| & \mbox{if } d_\Lambda=1,\\
|\mathbf{v}_1\times \mathbf{v}_2| & \mbox{if } d_\Lambda=2.
\end{cases}
$$
We let $\mathbf{w}_k$, $k=1, \cdots, d_\Lambda$, denote vectors in $\mathbb{R}^d$ satisfying
$$\mathbf{v}_j\cdot\mathbf{w}_k = 2\pi\delta_{jk},$$
that is, $\mathbf{w}_k$ are the primitive vectors of the reciprocal lattice $\Lambda^*$. Then the lattice vectors indexed by $m$ in the lattice $\Lambda$ and the reciprocal lattice $\Lambda^*$ are expressed by
\begin{align*}
\mathbf{R}_m:=\sum_{j=1}^{d_\Lambda} m_j\mathbf{v}_j \quad \mbox{and} \quad  \mathbf{K}_m:=\sum_{k=1}^{d_\Lambda} m_k\mathbf{w}_k \quad\mbox{for all }m=(m_1,\ldots,m_{d_\Lambda})\in\mathbb{Z}^{d_{\Lambda}},
\end{align*}
respectively. Let $\rho_m=\rho_m(\mathbf{x}):=|\mathbf{x}-\mathbf{R}_m|$, $\mathbf{\alpha}_m:=\mathbf{\alpha}+\mathbf{K}_m$, and $\beta_m=\beta_m(\kappa):=\sqrt{\kappa^2-|\mathbf{\alpha}_m|^2}$ with $\text{Im}({\beta}_m)\geq0$.

\subsection{The Ewald method} \label{subsec:approx:Ewald}
The Ewald method \cite{Ewald:1921:BOE} is known for efficient computation of the quasi-periodic Green function in free space (see \cite{Valerio:2007:CAA} for a comparative analysis with other computational methods). We examine the Ewald method to evaluate the quasi-periodic Green function $G^q_\kappa$, following the approach in \cite{Linton:2010:LSH}. For sufficiently large $n \in \N$, we approximate $G^q_\kappa$ by 
\begin{align} \label{def:qpG:E} 
    G^{q,E,n}_{\kappa}=G^{q,n}_{\kappa,\Lambda}+G^{q,n}_{\kappa,\Lambda^*}.
\end{align} 
For each $n \in \N$ and some splitting parameter $\eta>0$, the first term $G_{\kappa,\Lambda}^{q,n}$ is defined by
\begin{align*}
&G_{\kappa,\Lambda}^{q,n} (\mathbf{x}):=\frac{1}{4\pi}\sum_{|\mathbf{R}_m|\le n,\, m\in\mathbb{Z}^{d_\Lambda}}e^{{\rm i}\mathbf{\alpha}\cdot\mathbf{R}_m} \sum_{j=0}^\infty \frac{1}{j!} \left(\frac{\kappa}{2\eta} \right)^{2j} E_{j+1}\left(\eta^2 \rho_m^2 \right)
\end{align*}
for $d=2$, and by
\begin{align*}
&G_{\kappa,\Lambda}^{q,n} (\mathbf{x}):=\frac{1}{8\pi}\sum_{|\mathbf{R}_m|\le n,\, m\in\mathbb{Z}^{d_\Lambda}}\frac{e^{{\rm i}\mathbf{\alpha}\cdot\mathbf{R}_m}}{\rho_m}\left[e^{{\rm i}\kappa \rho_m}\operatorname{erfc}\left(\eta\rho_m+{\rm i}\frac{\kappa}{2\eta}\right)+e^{-{\rm i}\kappa \rho_m}\operatorname{erfc}\left(\eta\rho_m-{\rm i}\frac{\kappa}{2\eta}\right)\right]
\end{align*}
for $d=3$. Also, the second term $G_{\kappa,\Lambda^*}^{q,n}$ is defined by 
\begin{align*}
&G_{\kappa,\Lambda^*}^{q,n} (\mathbf{x}):=\frac{1}{4\pi A}\sum_{|\mathbf{K}_m|\le n,\, m\in\mathbb{Z}}e^{{\rm i}\mathbf{\alpha}_m\cdot\mathbf{x}}\sum_{j=0}^\infty \frac{\left(-\eta^2(x_2^2+x_3^2)\right)^j}{j!}E_{j+1}\left(-\frac{\beta_m^2}{4\eta^2} \right) 
\end{align*}
for $d=3$ and $d_\Lambda= 1$, and by
\begin{align*}
&G_{\kappa,\Lambda^*}^{q,n} (\mathbf{x}):=\frac{\rm i}{4A}\sum_{|\mathbf{K}_m|\le n,\, m\in\mathbb{Z}^{2}}\frac{e^{{\rm i}\mathbf{\alpha}_m\cdot\mathbf{x}}}{\mathbf{\beta}_m}\left[e^{{\rm i} \mathbf{\beta}_m x_d}\operatorname{erfc}\left(-\eta x_d-{\rm i}\frac{\mathbf{\beta}_m}{2\eta}\right)+e^{-{\rm i}\mathbf{\beta}_m x_d}\operatorname{erfc}\left(\eta x_d-{\rm i}\frac{\mathbf{\beta}_m}{2\eta}\right)\right]
\end{align*}
for $d=2,3$ and $d_\Lambda=d-1$. Here, $E_{j}$ is the exponential integral function of order $j$, and $\operatorname{erfc}$ is the complementary error function. We note that  $E_j$ with $j\geq 1$ satisfies 
$$
\left|E_j(z)\right| = \left| \int_z^\infty \frac{e^{-t}}{t^j} dt \right|
\leq \frac{Ce^{-\Re(z)}}{|\Re(z)^j|},\quad\forall z\in\mathbb{C}\backslash{\rm i}\mathbb{R},
$$
and $\operatorname{erfc}$ satisfies 
\begin{equation}\label{eq:bound:erfc}
    \left|\operatorname{erfc}(z)\right|=\frac{2}{\sqrt{\pi}}\left|\int_z^\infty e^{-t^2}\,dt\right|\le C e^{-\Re\{z^2\}},\quad\forall z\in\mathbb{C},
\end{equation}
where the bounds are obtained by decomposing $\int_{z}^{\infty}$ as $\int_{z}^{\Re z}+\int_{\Re z}^\infty$.
In the Ewald method, the splitting parameter $\eta$ determines the convergence rate of the two series $G_{\kappa,\Lambda}^{q,n}$ and $G_{\kappa,\Lambda^*}^{q,n}$.

We fix $D$ as a bounded subset of $\mathbb{R}^d$, choose an arbitrary $\kappa_0\in\mathbb{C}$, and set $K$ to be a bounded neighborhood of $\kappa_0$.
In the following, we estimate the decay rate of $|G_{\kappa}^q(\mathbf{x})-G_{\kappa}^{q,E,n}(\mathbf{x})|$ for $\mathbf{x}\in D$ and $\kappa\in K$ as $n\to\infty$.
For $d=3$, from triangle inequality and \eqref{eq:bound:erfc}, we have
\begin{align*}
\left|G_{\kappa,\Lambda}^{q,\infty} (\mathbf{x})-G_{\kappa,\Lambda}^{q,n} (\mathbf{x})\right|&\le \frac{1}{8\pi}\sum_{\left|\mathbf{R}_m\right|> n,\, m\in\mathbb{Z}^{d_\Lambda}}\frac{C}{\rho_m}\left[\left|\operatorname{erfc}\left(\eta\rho_m+{\rm i}\frac{\kappa}{2\eta}\right)\right|+\left|\operatorname{erfc}\left(\eta\rho_m-{\rm i}\frac{\kappa}{2\eta}\right)\right|\right]\\
&\le \frac{1}{8\pi}\sum_{\left|\mathbf{R}_m\right|> n,\, m\in\mathbb{Z}^{d_\Lambda}}\frac{2C}{\rho_m}e^{-\eta^2\rho_m^2+\frac{\kappa^2}{4\eta^2}}.
\end{align*}
Now fix any $\delta \in (0,1)$. Using the relation $\rho_m\ge |\mathbf{R}_m|-|\mathbf{x}|$, for sufficiently large $n$, we arrive at 
$$\left|G_{\kappa,\Lambda}^{q,\infty} (\mathbf{x})-G_{\kappa,\Lambda}^{q,n} (\mathbf{x})\right|\le \frac{C'(\delta)e^{\operatorname{Re}\kappa_0^2/(4\eta^2)}}{n e^{(1-\delta) n^2\eta^2}}\quad\mbox{for all }\mathbf{x}\in D\mbox{ and }\kappa\in K,$$
where $C'(\delta)$ is a constant that is independent of $\mathbf{x}$, $\kappa$ and $n$.
Applying the same argument for $d=2$, we obtain
$$\left|G_{\kappa,\Lambda}^{q,\infty} (\mathbf{x})-G_{\kappa,\Lambda}^{q,n} (\mathbf{x})\right|\le \frac{C'(\delta)e^{\operatorname{Re}\kappa_0^2/(4\eta^2)}}{\eta^2 \cdot n^2 e^{(1-\delta) n^2\eta^2}}\quad\mbox{for all }\mathbf{x}\in D\mbox{ and }\kappa\in K.$$
Similarly, since $\beta_m \approx {\rm{i}} |\mathbf{K}_m|$ 
for $|\mathbf{K}_m|>n$ with sufficiently large $n$, we have
$$\left|G_{\kappa,\Lambda^*}^{q,\infty} (\mathbf{x})-G_{\kappa,\Lambda^*}^{q,n} (\mathbf{x})\right|\le \frac{C''(\delta)\cdot 4\eta^2 }{n^2e^{(1-\delta)n^2/(4\eta^2)}}\quad\mbox{for all }\mathbf{x}\in D\mbox{ and }\kappa\in K$$
when $d=3$ and $d_\Lambda =1$, 
and $$\left|G_{\kappa,\Lambda^*}^{q,\infty} (\mathbf{x})-G_{\kappa,\Lambda^*}^{q,n} (\mathbf{x})\right|\le \frac{C''(\delta)}{n e^{{(1-\delta) n^2}/{(4\eta^2)}}}\quad\mbox{for all }\mathbf{x}\in D\mbox{ and }\kappa\in K$$
when $d=2,3$ and $d_\Lambda =d-1$.

Therefore, $G_{\kappa}^{q,E,n} = G_{\kappa,\Lambda}^{q,n}+G_{\kappa,\Lambda^*}^{q,n}$ converges uniformly to $ G_{\kappa}^q$ as $n\to\infty$ with, for large $n$,
\begin{align*}
\max_{\mathbf{x}\in D,\kappa\in K}\left|G_{\kappa}^q(\mathbf{x})-G_{\kappa}^{q,E,n}(\mathbf{x})\right|&\le \max_{\mathbf{x}\in D,\kappa\in K} \left|G_{\kappa,\Lambda}^{q,\infty} (\mathbf{x})-G_{\kappa,\Lambda}^{q,n} (\mathbf{x})\right| + \max_{\mathbf{x}\in D,\kappa\in K} \left|G_{\kappa,\Lambda^*}^{q,\infty} (\mathbf{x})-G_{\kappa,\Lambda^*}^{q,n} (\mathbf{x})\right|\\
&\le C(\delta)\left(e^{\operatorname{Re}\kappa_0^2/(4\eta^2)}e^{-(1-\delta)n^2\eta^2} + e^{-(1-\delta)n^{2}/(4\eta^2)}\right),
\end{align*}
where we $\delta\in(0,1)$ is arbitrary.
By a similar proof, we arrive at the following lemma.
\begin{lemma}\label{lemma:G:decayrate:Ewald}
For $i=0,1,2$ and $0<\delta<1$, there exists $N\in\mathbb{N}$ such that for all $n>N$,
\beq\label{eq:decayrate:Ewald}
\begin{aligned}
&\max_{\mathbf{x}\in D,\kappa\in K}\left|\nabla^i(G_{\kappa}^q(\mathbf{x})-G_{\kappa}^{q,E,n}(\mathbf{x}))\right|
\le C(\delta)\left(e^{\operatorname{Re}\kappa_0^2/(4\eta^2)}e^{-(1-\delta)n^2\eta^2} + e^{-(1-\delta)n^{2}/(4\eta^2)}\right).
\end{aligned}
\eeq
\end{lemma}

\subsection{The shifted Green function method} \label{subsec:approx:shifted}
Since the quasi-periodic Green function \eqref{eq:Gq} decays slowly at infinity and its defining series converges only conditionally, a half-space shifted Green function can be introduced to ensure the convergence with a better decay rate. 

Consider a shift parameter $h>0$, where $h \neq 2\pi/\beta_{\mathbf{0}}$ to avoid resonance \cite{Bruno:2014:RCT}, and a unit vector $\hat{\mathbf{v}}$ orthogonal to $\text{span}(\Lambda)$. The half-space shifted Green function of order $j \in \mathbb{Z}_{\geq 0}$ is given by
\begin{align} \label{eq:shiftedG}
G_{\kappa, j}(\mathbf{x}) = \sum_{l=0}^j (-1)^l  \binom{j}{l} G_{\kappa} (\mathbf{x}+lh\hat{\mathbf{v}}).
\end{align}
The corresponding shifted quasi-periodic Green function is defined by
\begin{align} \label{eq:shiftedGq}
G^q_{\kappa, j}(\mathbf{x}) = \sum_{m \in \mathbb{Z}^{d_\Lambda}} e^{\rm{i}\alpha \cdot \mathbf{R}_m} G_{\kappa,j} (\mathbf{x}-\mathbf{R}_m).
\end{align}
By \eqref{eq:shiftedG}, the quasi-periodic Green function $G^q_\kappa$ in \eqref{eq:Gq} is written as follows:
\begin{align}\label{eq:Gq:shifted} 
G_\kappa^{q}(\mathbf{x})   = G_{\kappa,j}^q(\mathbf{x})  - \sum_{l=1}^j (-1)^l  \binom{j}{l} G_{\kappa}^q (\mathbf{x}+lh\hat{\mathbf{v}}).
\end{align}
In particular, we have
\begin{align*} 
G_\kappa^q(\mathbf{x})  =  \sum_{m \in \mathbb{Z}^{d_\Lambda}} e^{\rm{i}\alpha \cdot \mathbf{R}_m} G_{\kappa,j} (\mathbf{x}-\mathbf{R}_m) -\frac{\rm i}{4A} \sum_{m \in \mathbb{Z}^{d_\Lambda}} \sum_{l=1}^j (-1)^l  \binom{j}{l} e^{\rm i \alpha_m \cdot \mathbf{x}} H_0^{(1)}(\beta_m |\mathbf{x}^\perp +lh\hat{v}|)
\end{align*}
for $d=3$ and $d_{\Lambda}=1$, and
\begin{align*} 
G_\kappa^q(\mathbf{x})  =  \sum_{m \in \mathbb{Z}^{d_\Lambda}} e^{\rm{i}\alpha \cdot \mathbf{R}_m} G_{\kappa,j} (\mathbf{x}-\mathbf{R}_m) -\frac{\rm i}{2A} \sum_{m \in \mathbb{Z}^{d_\Lambda}} \sum_{l=1}^j (-1)^l  \binom{j}{l} \frac{e^{{\rm{i}} \alpha_m \cdot \mathbf{x} + {\rm i} \beta_m \left|x_d+lh\right|}}{\beta_m},
\end{align*}
for $d=2,3$ and $d_{\Lambda}=d-1$, provided that $\mathbf{x}^\perp \neq -lh\hat{\mathbf{v}}$ for $l = 1,\ldots, j$ and that $(\kappa,\alpha,\Lambda)$ is not a Rayleigh--Wood anomaly triple \cite{Bruno:2017:RCQ, Bruno:2020:EQG}. 

Observe that, in the expression \eqref{eq:Gq:shifted}, the spatial singularities and the Rayleigh--Wood anomalies are separated into the first and the second terms, respectively. This representation yields rapid convergence with respect to the order $j$ while retaining all spatial and spectral singular terms \cite{Bruno:2014:RCT,Bruno:2017:TDQ,Bruno:2017:RCQ}.

In the shifted Green function method, for sufficiently large $n\in\N$, we approximate $G^q_\kappa$ by 
\begin{align} \label{def:qpG:S} 
    G_{\kappa}^{q,S,n}(\mathbf{x}) = G_{\kappa,j}^{q,n}(\mathbf{x}) + G_{\kappa,j}^{q,\rm{rem},n}(\mathbf{x})
\end{align}
The first term $G_{\kappa,j}^{q,n}$ is defined by
\begin{align*}
    G_{\kappa,j}^{q,n}(\mathbf{x}) = \sum_{|\mathbf{R}_m|\le n,\,m \in \mathbb{Z}^{d_\Lambda}} e^{\rm{i}\alpha \cdot \mathbf{R}_m} G_{\kappa,j} (\mathbf{x}-\mathbf{R}_m),
\end{align*}
and the second term $G_{\kappa,j}^{q,\rm{rem},n}$ is defined by 
\begin{align*}
G_{\kappa,j}^{q,\rm{rem},n}(\mathbf{x})  =  -\frac{\rm i}{4A} \sum_{|\mathbf{K}_m|\le n,\,m \in \mathbb{Z}} \sum_{l=1}^j (-1)^l  \binom{j}{l} e^{{\rm i} \alpha_m \cdot \mathbf{x}} H_0^{(1)}(\beta_m |\mathbf{x}^\perp + lh\hat{v}|)
\end{align*}
if $d=3$ and $d_\Lambda =1$, and
\begin{align*}
G_{\kappa,j}^{q,\rm{rem},n}(\mathbf{x})  = -\frac{\rm i}{2A}\sum_{|\mathbf{K}_m|\le n,\,m \in \mathbb{Z}^{d_\Lambda}}  \sum_{l=1}^j (-1)^l \binom{j}{l} \frac{e^{{\rm i} \alpha_m \cdot \mathbf{x} + {\rm i}  \beta_m |x_d+lh|}}{\beta_m}
\end{align*}
if $d=2,3$ and $d_{\Lambda}=d-1$.

Let $D$ be a bounded subset of $\RR^d$, choose any $M>0$ and $\kappa_0 \in \C$, and set $K$ to be a bounded neighborhood of $\kappa_0$. We now estimate the error between $G_\kappa^{q}(\mathbf{x})$ and $G_\kappa^{q,S,n}(\mathbf{x})$ for any $D_M:=\mathbf{x} \in D \cap \{\mathbf{x}:|\mathbf{x} \cdot \hat{\mathbf{v}}| \leq M\}$ and $\kappa \in K$ when $n$ is sufficiently large. We remark that, if $|\mathbf{x} \cdot \hat{\mathbf{v}}| \leq M$, there exists a constant $C(h,M)$ such that 
\[
\left|G_{\kappa,j}(\mathbf{x}) \right| \leq C(h,M) \big|\mathbf{x}^{\parallel} \big|^{-j^\flat}
\]
where $j^\flat$ is given by $(j+1)/2$ if $j$ is even and by $j/2+1$ if $j$ is odd; see \cite{Bruno:2014:RCT, Bruno:2017:TDQ} for $d=2$ and $d=3$, respectively. It follows that, as in subsection \ref{subsec:approx:Ewald}, for any $\delta \in (0,1)$ and large $n$,
\begin{align*}
&\left|G_\kappa^{q}(\mathbf{x}) - G_\kappa^{q,S,n}(\mathbf{x})\right| \\ \qquad
& \leq
\sum_{\left|\mathbf{R}_m\right|> n,\,m \in \mathbb{Z}^{d_\Lambda}}  C(h,M) \big|(\mathbf{x}-\mathbf{R}_m)^\parallel\big|^{-j^\flat} + C'\sum_{\left|\mathbf{K}_m\right|> n,\,m \in \mathbb{Z}^{d_\Lambda}}  2^j \cdot \frac{e^{-\text{Im}(\beta_m) |\mathbf{x}^\perp|}}{\left|\beta_m\right|^{d_\Lambda/2}}\\
& \leq C(h,M,j) n^{-j^\flat +d_\Lambda+1} + C'(\delta,j)n^{-d_\Lambda/2} e^{-(1-\delta)n |\mathbf{x}^\perp|}
\end{align*} 
provided that $j^\flat > d_\Lambda+1$ and that, when $d_\Lambda =1$, $x^\perp$ bounded away from $0$. Hence, $G_{\kappa}^{q,S,n}$ converges uniformly to $ G_{\kappa}^q$ as $n\to\infty$ when the order $j$ is larger than $2d_\Lambda$. Ias in Lemma \ref{lemma:G:decayrate:Ewald}. As in Lemma \ref{lemma:G:decayrate:Ewald}, we can similarly establish the estimate below for the shifted Green function method. We assume that there exists a constant \(b_D\) such that \(|\mathbf{x}^\perp|\ge b_D\) for all \(\mathbf{x}\in D_M\), where $b_D>0$ if \(d_\Lambda=1\),
\begin{lemma}\label{lemma:G:decayrate:shifted}
For $i=0,1,2$ and $0<\delta<1$, there exists $N\in\mathbb{N}$ such that for all $n>N$,
\beq\label{eq:decayrate:shifted}
\begin{aligned}
\max_{\mathbf{x}\in D_M,\kappa\in K}\left|\nabla^i(G_\kappa^{q}(\mathbf{x}) - G_\kappa^{q,S,n}(\mathbf{x}))\right|\le C(h,M,j) n^{-j^\flat +d_\Lambda+1} + C'(\delta,j)n^{-d_\Lambda/2} e^{-(1-\delta)b_Dn}
\end{aligned}
\eeq
provided that $j > 2d_\Lambda$.
\end{lemma}
\noindent This follows from   
$\left|\nabla^i G_{\kappa,j}(\mathbf{x}) \right| \leq C_i(h,M) \big|\mathbf{x}^{\parallel}\big|^{-j^\flat-i}$,
which can be derived from Section 3(a) of Bruno et al. \cite{Bruno:2017:TDQ}.

\section{Approximation of a resonant frequency} \label{sec:resonant}

\subsection{Approximation schemes}
Based on the approximation methods for the quasi-periodic Green function discussed in the previous section, we develop approximation schemes for a resonant frequency of the problem \eqref{eq:eigproblem:simple}. We recall that the resonant frequency corresponds to the characteristic value of the operator-valued function $\mathcal{A}(\omega)$ defined in \eqref{def:A}.

We use {\rm M} to denote the approximation method for quasi-periodic Green function, with ${\rm M}={\rm E}$ for the Ewald method and ${\rm M}={\rm S}$ for the shifted Green function method.
Let $\omega_0\in\mathbb{C}$ be a characteristic value of $\mathcal{A}(\omega)$.
By Lemma \ref{lemma:poles:characterization}, we have $\omega_0\in\mathbb{R}$ and there exists a $\eps>0$ such that $\omega=\omega_0$ is the only characteristic value of $\mathcal{A}(\omega)$ in $B_\eps(\omega_0)$.
We define, for ${\rm M}\in\{{\rm E},{\rm S}\}$,    
\begin{align}  \label{def:An} 
\mathcal{A}_n^{\rm M}(\omega):=\begin{pmatrix}
\mathcal{S}_{D}^{n,{\rm M}}(n_1\omega) & -\mathcal{S}_{D}^{n,{\rm M}}(n_0\omega)
\\
C_{01}^2\left(-\frac{1}{2}{\rm id}+\mathcal{K}_{D}^{n,{\rm M},*}(n_1\omega)\right)&-\frac{1}{2}{\rm id}-\mathcal{K}_{D}^{n,{\rm M},*}(n_0\omega)
\end{pmatrix},
\end{align}
where $\mathcal{S}_{D}^{n,{\rm M}}$ and $\mathcal{K}_{D}^{n,{\rm M},*}$ are defined by the right-hand sides of \eqref{def:S} and \eqref{def:Kstar} with the Green function $G_\kappa^q$ replaced by $G_{\kappa,\Lambda}^{q,n}+G_{\kappa,\Lambda^*}^{q,n}$ for the Ewald method (${\rm M}={\rm E}$) and $G_{\kappa,j}^{q,n}+G_{\kappa,j}^{q,{\rm rem},n}$ for the shifted Green function method (${\rm M}={\rm S}$).

\begin{thm}\label{thm:convergence}
Let $\mu$ be the multiplicity of $\omega_0$. There exists an $N\in\mathbb{N}$ such that for all $n>N$, each of the operator-valued functions $\mathcal{A}_n^{\rm E}(\omega)$ and $\mathcal{A}_n^{\rm S}(\omega)$ has $\mu$ characteristic values in $B_\eps(\omega_0)$ counting multiplicity, namely, $\{\omega_{n,k}^{\rm E}\}_{k=1}^\mu$ and $\{\omega_{n,k}^{\rm S}\}_{k=1}^\mu$. Moreover, there exists a constant $C>0$ satisfying, for all $n>N$,
\begin{align}
|\omega_{n,j}^{\rm E}-\omega_0|&\le C \left(e^{\operatorname{Re}\kappa_0^2/(4\eta^2)}e^{-(1-\delta)n^2\eta^2} + e^{-(1-\delta)n^{2}/(4\eta^2)}\right),\\
|\omega_{n,j}^{\rm S}-\omega_0|&\le C \left(n^{-j^\flat +d_\Lambda+1} + n^{-d_\Lambda/2} e^{-(1-\delta)b_D n}\right).
\end{align}
\end{thm}

\subsection{Proof of convergence}  \label{subsec:resonant:proof}
In this subsection, we provide a proof of Theorem \ref{thm:convergence}. The key idea behind the proof is to leverage Gohberg--Sigal theory \cite{Gohberg:1971:OGL}. For applications of Gohberg--Sigal theory with layer potential techniques, see \cite{Ammari:2004:SRS, Ammari:2009:LPT}. We also need the following lemma.

\begin{lemma}\label{lemma:supnormconv}
Let $O$ be a bounded subset of $\mathbb{C}$. For both ${\rm M}\in\{{\rm E},{\rm S}\}$, there holds $$\textstyle\lim_{n\to\infty}\sup_{\omega\in O}\|\mathcal{A}(\omega)-\mathcal{A}_n^{\rm M}(\omega)\|_{\mathcal{L}(L^2,L^2)}=0.$$
\end{lemma}
\begin{proof}
For both ${\rm M}\in\{{\rm E},{\rm S}\}$ and all $n\in\mathbb{N}$, we have
$$\mathcal{A}(\omega)-\mathcal{A}_n^{\rm M}(\omega) = \begin{pmatrix}
\mathcal{S}_{D}^{q}(n_1\omega)-\mathcal{S}_{D}^{n,{\rm M}}(n_1\omega) & -\mathcal{S}_{D}^{q}(n_0\omega)+\mathcal{S}_{D}^{n,{\rm M}}(n_0\omega)
\\
C_{01}^2\left(\mathcal{K}_{D}^{q,*}(n_1\omega)-\mathcal{K}_{D}^{n,{\rm M},*}(n_1\omega)\right)&-\mathcal{K}_{D}^{q,*}(n_0\omega)+\mathcal{K}_{D}^{n,{\rm M},*}(n_0\omega)
\end{pmatrix}.$$
We fix $0<\delta<1$ arbitrarily.
From Lemma \ref{lemma:G:decayrate:Ewald}, for $n$ large enough, we have
$$\sup_{\omega\in O}\|\mathcal{A}(\omega)-\mathcal{A}_n^{\rm E}(\omega)\|_{\mathcal{L}(L^2,L^2)}\le C(\delta)\left(e^{\operatorname{Re}\kappa_0^2/(4\eta^2)}e^{-(1-\delta)n^2\eta^2} + e^{-(1-\delta)n^2/(4\eta^2)}\right).$$
From Lemma \ref{lemma:G:decayrate:shifted}, for $n$ large enough, we have
$$\sup_{\omega\in O}\|\mathcal{A}(\omega)-\mathcal{A}_n^{\rm S}(\omega)\|_{\mathcal{L}(L^2,L^2)}\le C(h,M,j) n^{-j^\flat +d_\Lambda+1} + C'(\delta,j)n^{-d_\Lambda/2} e^{-(1-\delta)b_D n}.$$
The right-hand sides converge to zero as $n\to\infty$.
\end{proof}

We now prove the main theorem. 

\medskip
\noindent{\bf Proof of Theorem \ref{thm:convergence}.} 
For either $M \in \{{\rm E}, {\rm S}\}$, we write $\mathcal{A}_n=\mathcal{A}_n^{\rm M}$ for simplicity.
$\mathcal{A}(\omega)$ is invertible for all $\omega\in\p B_\eps(\omega_0)$.
From Lemma \ref{lemma:supnormconv}, we may choose $N$ large enough so that
\begin{equation}\label{eq:def:rho}
\rho_n:=\sup_{\omega\in\p B_\eps(\omega_0)}\|\mathcal{A}(\omega)^{-1}(\mathcal{A}(\omega)-\mathcal{A}_n(\omega))\|<1/2\quad\mbox{for all }n>N.
\end{equation}
Thus, by Theorem \ref{thm:GS:Rouche} (the generalized Rouch\'{e} theorem), for all $n>N$, there exist $\mu$ characteristic values of $\mathcal{A}_n(\omega)$ in $B_\eps(\omega_0)$ counting multiplicity, namely, $\{\omega_{n,k}\}_{k=1}^\mu$.
From the generalized argument principle \cite[Theorem 2.2 with $f(z)=z^m$]{Beceanu:2026:EPL} (see also Theorem \ref{thm:GS:argument}), we have, for $m=1,\dots,\mu$, 
\begin{align*} 
\sum_{k=1}^\mu \left(\omega_{n,k}-\omega_0\right)^m 
= &\frac{1}{2\pi{\rm i}}\operatorname{tr}\int_{\p B_\eps(\omega_0)}\left(\omega-\omega_0\right)^m \mathcal{A}_n(\omega)^{-1}\frac{\p}{\p\omega}\mathcal{A}_n(\omega)\,d\omega\\ 
=&\frac{1}{2\pi{\rm i}}\sum_{p=0}^\infty \operatorname{tr}\int_{\p B_\eps(\omega_0)}\left(\omega-\omega_0\right)^m \left[\mathcal{A}(\omega)^{-1}(\mathcal{A}(\omega)-\mathcal{A}_n(\omega))\right]^p \mathcal{A}(\omega)^{-1}\frac{\p}{\p\omega}\mathcal{A}_n(\omega)\,d\omega.
\end{align*}
We then prove the convergence rate of $\omega_{n,j}\to\omega_0$ as $n\to\infty$.
Proceeding as in the proof of \cite[Theorem 3.9]{Ammari:2009:LPT}, one finds that for $m=1,\dots,\mu$,
\beq\label{eq:multifreq:argprin:2} 
\sum_{k=1}^\mu \left(\omega_{n,k}-\omega_0\right)^m  = \frac{1}{2\pi{\rm i}}\sum_{p=1}^\infty\frac{m}{p}\operatorname{tr}\int_{\p B_\eps(\omega_0)}\left(\omega-\omega_0\right)^{m-1}\left[\mathcal{A}(\omega)^{-1}(\mathcal{A}(\omega)-\mathcal{A}_n(\omega))\right]^p\,d\omega.
\eeq
Using the analogue of \cite[Lemma 4.2]{Beceanu:2026:EPL} and the facts that $\sum_{p=m}^\infty p x^p=\frac{m-(m-1)x}{(1-x)^2}x^m$ for $0\leq x < 1$ and that $m\le 2^{m-1}$ for all $1\le m\le\mu$ with \eqref{eq:def:rho}, we arrive at 
\begin{align}\notag
\left|\sum_{k=1}^\mu \left(\omega_{n,k}-\omega_0\right)^m \right| 
&\le\sum_{p=m}^\infty \frac{m}{p}p^2\mu\delta^m{\rho_{n}}^p=m\eps^m \mu\left(\frac{m-(m-1){\rho_{n}}}{(1-\rho_{n})^2}\right){\rho_{n}}^{m} \\ \notag
& \le  2^{m-1}\eps^m \mu^m\left(\frac{2^{m-1}}{(1-\rho_{n})^{2}}\right){\rho_{n}}^{m} \le\left(4\eps\mu\rho_{n}\right)^m\quad\mbox{for }m=1,\dots,\mu.
\end{align}
Finally, using \cite[Lemma 4.3]{Beceanu:2026:EPL}, we have
$$ \left|\omega_{n,k}-\omega_0\right|\le 8\eps\mu\rho_n\quad\mbox{for } k=1, \dots, \mu.
$$
Hence, the conclusion follows from \eqref{eq:def:rho} with Lemma \ref{lemma:G:decayrate:Ewald} and Lemma \ref{lemma:G:decayrate:shifted}.



\section{Numerical computation of quasi-periodic Green function and resonant frequencies} \label{sec:num} 
In this section, we demonstrate the numerical computation of resonance frequencies. We perform computations in two dimensions via the two proposed schemes, using the Ewald method and the shifted Green function method. Previous studies have considered the computation of resonance frequencies \cite{Ammari:2017:MNF, Ammari:2018:MRA} and methods for evaluating quasi-periodic Green functions \cite{Bruno:2020:EQG}. 
Here, using the Ewald and shifted Green function methods, we construct the corresponding operators $\mathcal{A}_n$ as described in Section \ref{sec:resonant} and obtain numerical approximations of the target resonance frequencies.

\begin{figure}[h] 
    \centering
    {\includegraphics[width=0.7\textwidth]{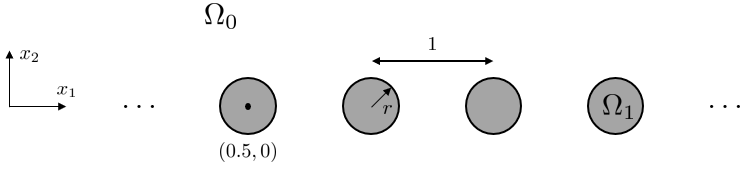}} 
    \caption{\label{fig:config}
	Configuration of the numerical simulations wth $r = 0.05$ or $0.1$.}
\end{figure}

For the two-dimensional case $d=2$, we have $d_\Lambda =1$. We set the lattice $\Lambda$ to be $\Z \times \{0\}$ so that the primitive vector is $\mathbf{v}_1 = (1,0)$, and define $Y = (0,1) \times \mathbb{R}$ as the primitive cell.
In all simulations, we take $D$ as a disk of radius $r = 0.05$ or $0.1$ centered at $(0.5,0) \in \mathbb{R}^2$. Given the refractive index $n_0$ and the angular frequency $\omega$,  we set $\alpha = (0.8 |n_0\omega|, 0) \in \mathbb{R}^2$ when designing the incident field (see Section \ref{subsec:prob:periodic}). 
The configuration of the simulation is described in Figure \ref{fig:config}.

In the following subsections, we presen computational examples and convergence behaviors of the quasi-periodic Green functions (Section \ref{subsec:num:qpGreen}) and resonant frequencies (Section \ref{subsec:num:resonant}).
All computations in this section are performed in MATLAB, with the boundary element method implemented using GYBSILAB \cite{Alouges:2018:FBS}. The boundaries were discretized using $N$ uniformly spaced vertices along each circumference, where $N=50$ for $r=0.05$ and $N=100$ for $r=0.1$. 

\subsection{Quasi-periodic Green function by Ewald and shifted Green function method}  \label{subsec:num:qpGreen} 
Recall that $\rm E$ and $\rm S$ denote the Ewald and shifted Green function methods. To compute quasi-periodic Green functions used in the construction of $\mathcal{A}_n^{\rm M}$, we employ the Ewald method (${\rm M} = {\rm E}$) and the shifted Green function method (${\rm M} = {\rm S}$) described in Sections\ref{subsec:approx:Ewald} and \ref{subsec:approx:shifted}. 

In the Ewald method \eqref{def:qpG:E}, we truncate the inner infinite sum in $G^{q,n}_{\kappa,\Lambda}$ as $J=10$ and compute it as follows:
\begin{align*}
\frac{1}{4\pi}\sum_{|\mathbf{R}_m|\le n,\, m\in\mathbb{Z}^{d_\Lambda}}e^{{\rm i}\mathbf{\alpha}\cdot\mathbf{R}_m} \sum_{j=0}^J \frac{1}{j!} \left(\frac{\kappa}{2\eta} \right)^{2j} E_{j+1}\left(\eta^2 \rho_m^2 \right).
\end{align*}
Here, the exponential integral functions $E_j$ are evaluated using the genexpint package \cite{Bocquet:2015:GEI}.
To compute $G^{q,n}_{\kappa,\Lambda^*}$, we use the implementation of the $\operatorname{erfc}$ function provided by the Faddeeva package \cite{Johnson:2012:FPC}, which supports complex-valued arguments. 

On the other hand, in the shifted Green function method \eqref{def:qpG:S}, we use the half-space shifted Green function $G_{\kappa,j}$ with order $j=8$; see \eqref{eq:shiftedG} for its explicit definition. The order $j$ is retained throughout the expression for the approximation $\eqref{def:qpG:S}$.

\begin{figure}[h!]
    \centering
    \subfloat[\centering Ewald method (log-linear scale)]{\includegraphics[width=.42\linewidth]{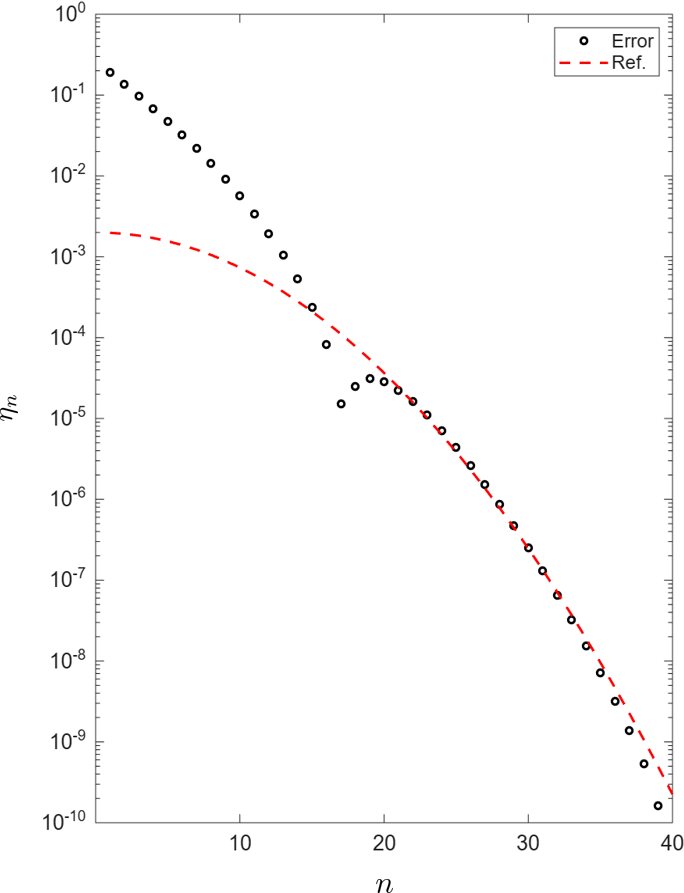}} 
    \qquad
    \subfloat[\centering Shifted Green function method (log-log scale)]{\includegraphics[width=.42\linewidth]{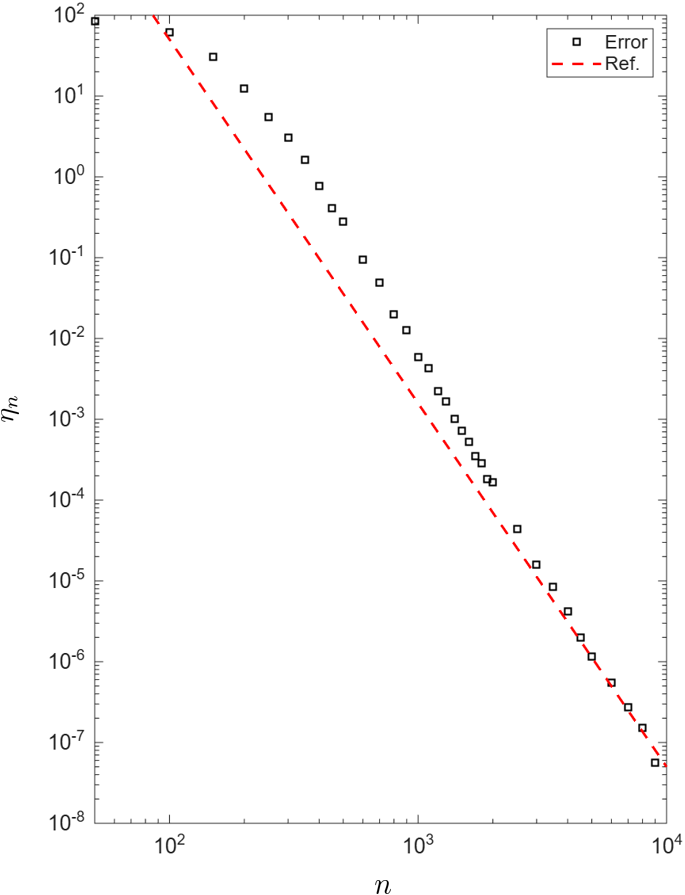}}
    \caption{\label{fig:qpGreen:xy}
    Relative errors $\eta_n(\mathbf{x},\mathbf{y})$ of the quasi-periodic Green function for $\mathbf{x} =(0.5,0)$ and $\mathbf{y} = (0.55,0.05)$, with $r=0.05$ (see \eqref{def:err:eta}). The errors obtained using the Ewald method are shown in (a) for $n=1,2,\dots,39$ on a log-linear scale, and those obtained using the shifted Green function method are shown in (b) for $n=50,100,\ldots,500$, $600, 700, \ldots, 2000$, $2500,3000,\ldots,5000$, $6000,7000,8000,9000$ on a log-log scale. The reference curves are $y = 0.02e^{-4n^2}$ in (a) and $y=5\cdot10^{10} n^{-4.5}$ in (b).} 
\end{figure}

\begin{figure}[h!]
    \centering
    \subfloat[\centering Ewald method(log-linear)]{\includegraphics[width=.42\linewidth]{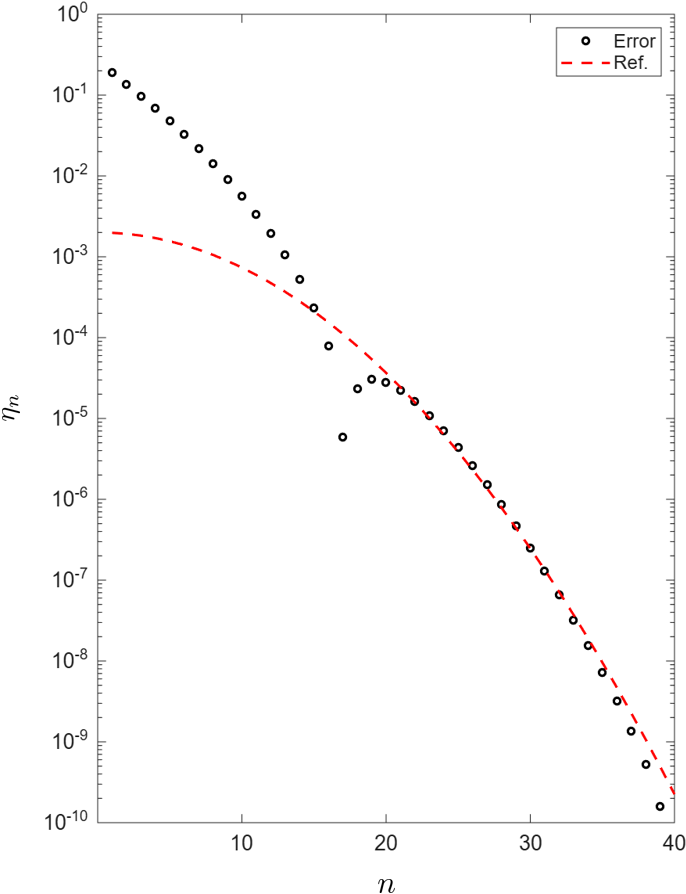}}
    \qquad
    \subfloat[\centering Shifted Green function method (log-log scale)]{\includegraphics[width=.42\linewidth]{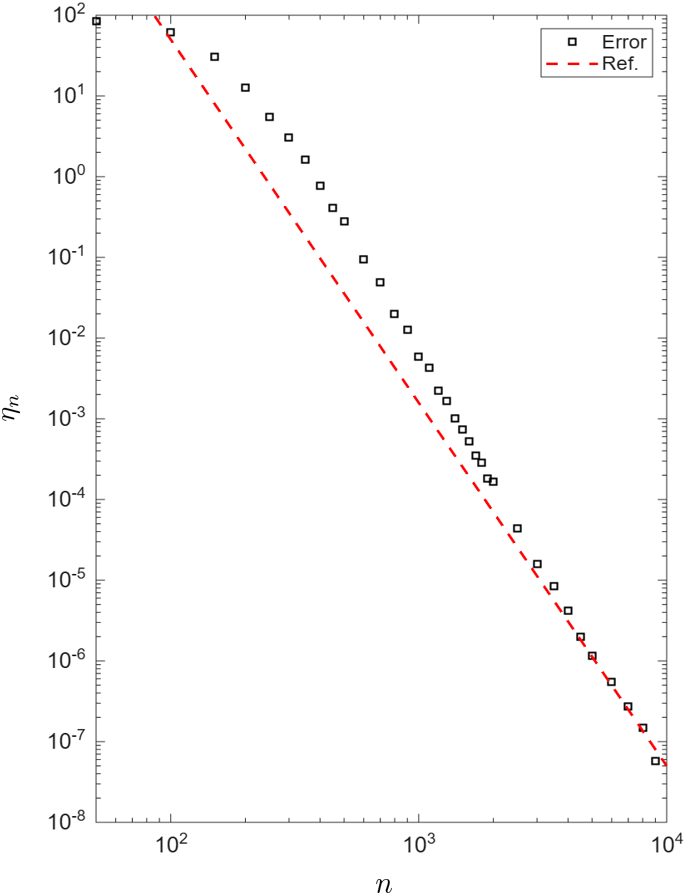}}
    \caption{\label{fig:qpGreen:xx}
    Relative errors $\eta_n(\mathbf{x},\mathbf{x})$ of the quasi-periodic Green function for $\mathbf{x} =(1/2,0)$, with $r=0.05$.The errors obtained using the Ewald method are shown in (a) for $n=1,2,\dots,39$ on a log-linear scale, and those obtained using the shifted Green function method are shown in (b) for $n=50,100,\ldots,500$, $600, 700, \ldots, 2000$, $2500,3000,\ldots,5000$, $6000,7000,8000,9000$ on a log-log scale. The reference curves are $y = 0.02e^{-4n^2}$ in (a) and $y=5\cdot10^{10} n^{-4.5}$ in (b).} 
\end{figure}

To observe the convergence behaviors of the quasi-periodic Green functions using the Ewald and shifted Green function method, we present the relative errors $\eta_n$ yielded by both methods in Figures \ref{fig:qpGreen:xy} and $\ref{fig:qpGreen:xx}$. Here, $\eta_n$ is given by
\begin{align} \label{def:err:eta}
    \eta_n(\mathbf{x} ,\mathbf{y}) := \left|\frac{G^{q,{\rm M},n}_\kappa(\mathbf{x}-\mathbf{y}) - G^{q,{\rm M},n_\infty}_\kappa(\mathbf{x}-\mathbf{y}) }{G^{q,{\rm M},n_\infty}_\kappa(\mathbf{x}-\mathbf{y})} \right|,
\end{align}
where $n_\infty = 40$ for $\rm{M} = {\rm E}$, and $n_\infty = 10^4$ for $\rm{M} = {\rm S}$ are used as the references.
For the points ${\bf x} = (1/2,0)$ and ${\bf y} = (0,1)$, Figures \ref{fig:qpGreen:xy} and \ref{fig:qpGreen:xx} show the relative errors $\eta_n(\bf{x},\bf{y})$ and $\eta_n(\bf{x},\bf{x})$, where $n_0\omega=0.05$ and $n_1\omega = 0.1$. The Ewald method (Figures \ref{fig:qpGreen:xy}(a) and \ref{fig:qpGreen:xx}(a)) are shown in log-linear scale, and the shifted Green function methods (Figures \ref{fig:qpGreen:xy}(b) and \ref{fig:qpGreen:xx}(b)) are shown in log-log scale. We remark that the errors decrease exponentially by the Ewald method, while algebraically by the shifted Green function method, which is consistent with the convergence results of Section \ref{sec:Green}.

\subsection{Resonant frequencies} \label{subsec:num:resonant} 
\subsubsection{Computation scheme}

We aim to find the resonant frequency $\omega_0$ for which the operator $\mathcal{A}(\omega_0)$ becomes singular. Accordingly, we approximate $\mathcal{A}(\omega)$ by $\mathcal{A}_n$ for sufficiently large $n$, and investigate the smallest modulus of the eigenvalues
\begin{align} \label{def:lambda0}
    \lambda_0(\omega) := \min_{\lambda \in \sigma(\mathcal{A}_n(\omega))}|\lambda|,
\end{align}
where $\sigma(\mathcal{A}_n(\omega))$ denotes the spectrum of $\mathcal{A}_n(\omega)$. Since Lemma \ref{lemma:poles:characterization} establishes that the resonant frequencies are real-valued, we evaluate $\mathcal{A}_n(\omega)$ only for real frequencies $\omega$.

The process of identifying the resonant frequency is as follows. 
\begin{itemize}
    \item {\bf Step 1.} Choose a sufficiently large truncation size $n=n_k$ by examining the relative error of the quasi-periodic Green functions at some sample points of $\kappa, \mathbf{x}, \mathbf{y}$. 
    Specifically, we choose $k$ such that the following criterion is satisfied:
    \begin{align} 
        \left|\frac{G^{q,{\rm M},n_{k+1}}_\kappa(\mathbf{x}-\mathbf{y}) - G^{q,{\rm M},n_k}_\kappa(\mathbf{x}-\mathbf{y}) }{G^{q,{\rm M},n_{k+1}}_\kappa(\mathbf{x}-\mathbf{y})} \right| < 10^{-7},
    \end{align}
    where $n_k = 5k$ for $\rm{M} = {\rm E}$ and $n_k= 100k$ for $\rm{M} = {\rm S}$. 
    
    \item {\bf Step 2.} Compute the eigenvalues of $\mathcal{A}_n(\omega)$ for some candidate $\omega$ on the real line, and find $\omega_0$ that yields the minimal modulus (see \eqref{def:lambda0}) using the bisection method. In this step, we determine $\omega_0$ up to four decimal places.   
\end{itemize}

\subsubsection{Examples} \label{subsec:num:result}
Under the setting described in the previous subsection, we identify the resonant frequency $\omega_0$ for several pairs of refractive indices $n_0$ and $n_1$. 

\begin{figure}[] 
    \centering
    {\includegraphics[width=0.8\textwidth]{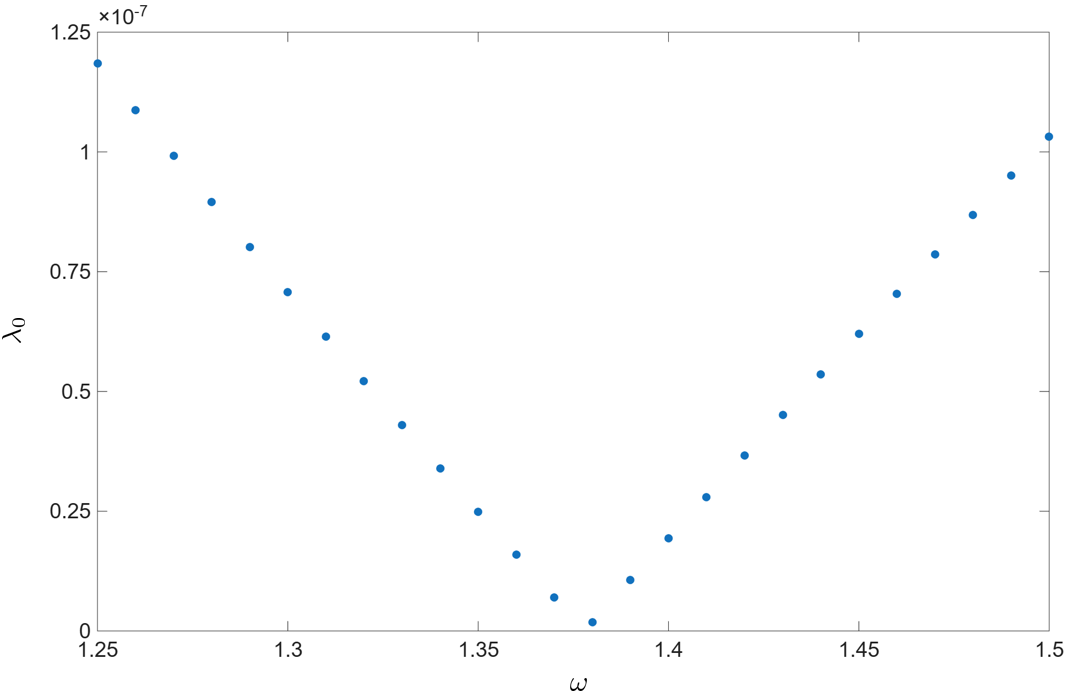}} 
    \caption{\label{fig:search}
	Minimum eigenvalue modulus $\lambda_0 =\lambda_0(\omega)$ (see \eqref{def:lambda0}) for $(n_0,n_1) = (0.005, 0.01)$, with $r=0.05$, ${\rm M} = {\rm E}$, and $n=40$.}
\end{figure}

For each pair of refractive indices $(n_0,n_1) = (0.005,0.01), (0.005,0.1), (0.05,0.01), (0.05,0.1)$, we find a resonant frequency by examining the spectral quantity $\lambda_0(\omega)$ \eqref{def:lambda0} as $\omega$ varies over  $\mathbb{R}_{>0}$. The case $(n_0,n_1) = (0.005,0.01)$ is shown in Figure \ref{fig:search}, where $\omega$ varies over $[1.25, 1.5]$ with a step size $0.01$, and the minimum is attained at $\omega = 1.38$. Table \ref{tab:res} shows computed resonant frequencies $\omega_0$ up to four decimal places for the considered pairs of refractive indices. 

\begin{table}
\centering
\caption{\label{tab:res} 
Computed resonant frequencies for pairs of refractive indices, up to four decimal places, with the corresponding values of $\lambda_0$ \eqref{def:lambda0}, for $r=0.05$.}
\begin{tabular}{c | c c}
\hline
\diagbox{$n_0$}{$n_1$} & $0.01$ & $0.1$ \\ \hline  
$0.005$ & $1.3779$ ($\lambda_0 = 8.28 \times 10^{-12}$)  & $0.0930$ ($\lambda_0 = 2.56 \times 10^{-10}$)\\    
$0.05$ & $0.4841$ ($\lambda_0 = 1.38 \times 10^{-10}$) & $0.2911$ ($\lambda_0 = 6.40 \times 10^{-11}$)\\    
\hline
\end{tabular}
\end{table}

Lastly, we examine the convergence behavior of $\mathcal{A}_n(\omega_0)$ with respect to the truncation size $n$. In particular, with a suitable step size $\delta \in \mathbb{N}$, we estimate the relative errors  
\begin{align} \label{def:err:op}
    \frac{\left\|\mathcal{A}_{n-k\delta}(\omega_0) - \mathcal{A}_{n+\delta}(\omega_0)\right\|_2}{\left\|\mathcal{A}_{n+\delta}(\omega_0)\right\|_2}
\end{align}
for $k = 0,1,2,\dots$, where the truncation order $n+\delta$ is used as the reference and $\|\cdot\|_2$ denotes the operator norm induced by the $L^2$-norm. Here, we set $\delta = 5$ for $\rm{M} = {\rm E}$, and $\delta = 500$ for $\rm{M} = {\rm S}$.
Tables \ref{tab:err:op1}, \ref{tab:err:op2} and \ref{tab:err:op2:large} present the relative errors \eqref{def:err:op} of the operator $\mathcal{A}_n(\omega_0)$ for each pair of refractive indices, with radii $r=0.05$ and $r=0.1$. Across all three tables, the convergence behaviors of $\mathcal{A}^{\rm E}_n(\omega_0)$ and $\mathcal{A}^{\rm S}_n(\omega_0)$ differ, depending on the approximation method used for the quasi-periodic Green function (see Lemma \ref{lemma:supnormconv}). We also observe slower convergence for the larger radius ($r=0.1$; Table \ref{tab:err:op2:large}) than for the smaller radius ($r=0.05$; Table \ref{tab:err:op2}) for the same pair of refractive indices.

\begin{table}
\centering 
\caption{\label{tab:err:op1} 
Relative errors of $\mathcal{A}_n(\omega)$ (see \eqref{def:err:op}) for $(n_0,n_1) = (0.005,0.01)$, $\omega = 1.3779$, and $r=0.05$, computed based on the Ewald method with truncation order $n=55$ in (a) and the shifted Green function method with truncation order $n=12000$ in (b) as the respective references.}
\begin{subtable}{0.35\textwidth} 
    \centering
    \caption{Ewlad method} \vspace{0.5em}
    \begin{tabular}{c c}
        \hline
        $n$ & relative error \\ \hline
        25 & 0.0987\\ 
        30 & 0.0761\\
        35 & 0.0565\\ 
        40 & 0.0396\\
        45 & 0.0247\\ 
        50 & 0.0116\\
        \hline
    \end{tabular} 
\end{subtable}
\quad
\begin{subtable}{0.35\textwidth} 
    \centering
    \caption{Shifted Green function method} \vspace{0.5em}
    \begin{tabular}{c c}
        \hline
        $n$ & relative error \\ \hline
        9000 & 0.0979\\  
        9500 & 0.0748\\
        10000 & 0.0553\\ 
        10500 & 0.0385\\
        11000 & 0.0239\\ 
        11500 & 0.0111\\
        \hline
    \end{tabular} 
\end{subtable}
\end{table}

\begin{table}[h]
\centering 
\caption{\label{tab:err:op2} 
Relative errors of $\mathcal{A}_n(\omega)$ (see \eqref{def:err:op}) for $(n_0,n_1) = (0.05,0.1)$, $\omega = 0.2911$, and $r=0.05$, computed based on the Ewald method with truncation order $n=45$ in (a) and the shifted Green function method with truncation order $n=7000$ in (b) as the respective references.}
\begin{subtable}{0.35\textwidth} 
    \centering
    \caption{Ewlad method} \vspace{0.5em}
    \begin{tabular}{c c}
        \hline
        $n$ & relative error \\ \hline
        15 & 0.1093\\ 
        20 & 0.0780\\
        25 & 0.0540\\ 
        30 & 0.0353\\
        35 & 0.0206\\ 
        40 & 0.0090\\
        \hline
    \end{tabular} 
\end{subtable}
\quad
\begin{subtable}{0.35\textwidth} 
    \centering
    \caption{Shifted Green function method} \vspace{0.5em}
    \begin{tabular}{c c}
        \hline
        $n$ & relative error \\ \hline
        4000 & 0.1154\\ 
        4500 & 0.0522\\
        5000 & 0.0279\\ 
        5500 & 0.0226\\
        6000 & 0.0171\\ 
        6500 & 0.0086\\
        \hline
    \end{tabular} 
\end{subtable}
\end{table}

\begin{table}[h]
\centering 
\caption{\label{tab:err:op2:large} 
Relative errors of $\mathcal{A}_n(\omega)$ (see \eqref{def:err:op}) for $(n_0,n_1) = (0.05,0.1)$, $\omega = 0.2789$, and $r=0.1$, computed based on the Ewald method with truncation order $n=45$ in (a) and the shifted Green function method with truncation order $n=7000$ in (b) as the respective references.}
\begin{subtable}{0.35\textwidth} 
    \centering
    \caption{Ewlad method} \vspace{0.5em}
    \begin{tabular}{c c}
        \hline
        $n$ & relative error \\ \hline
        15 & 0.1109\\ 
        20 & 0.0797\\
        25 & 0.0555\\ 
        30 & 0.0364\\
        35 & 0.0210\\ 
        40 & 0.0095\\
        \hline
    \end{tabular} 
\end{subtable}
\quad
\begin{subtable}{0.35\textwidth} 
    \centering
    \caption{Shifted Green function method} \vspace{0.5em}
    \begin{tabular}{c c}
        \hline
        $n$ & relative error \\ \hline
        4000 & 0.1262\\ 
        4500 & 0.0726\\
        5000 & 0.0424\\ 
        5500 & 0.0290\\
        6000 & 0.0202\\ 
        6500 & 0.0100\\
        \hline
    \end{tabular} 
\end{subtable}
\end{table}

\section{Conclusion} \label{sec:conclusion}
We have investigated convergence analysis of resonant frequencies for periodic Helmholtz problems. By formulating the Helmholtz problem in a periodic domain as a system of boundary integral equations, the resonant frequencies correspond to the characteristic values of the operator-valued function in the equation. Using the Gohberg-Sigal theory, we proposed two approximation schemes for the resonant frequencies, whose convergence rates are governed by the underlying approximation of the quasi-periodic Green function. We also examined the convergence behaviors of the proposed schemes via numerical simulations. As a future direction, we plan to investigate numerical schemes for evaluating  the operator $\mathcal{A}_n(\omega)$ and examine their impact on the associated convergence behavior of the resonance frequencies.


\appendix
\section{Gohberg--Sigal theory on operator-valued functions} \label{appen:GS}
Gohberg--Sigal theory provides a generalization of the argument principle and Rouch\'e's theorem in complex analysis to operator-valued functions. For Banach spaces $\mathcal{B}$ and $\mathcal{B}'$, denote by $\mathcal{L}(\mathcal{B},\mathcal{B}')$ the algebra of all bounded linear operators from $\mathcal{B}$ to $\mathcal{B}'$. Let $\mathfrak{U}(z_0)$ denote the class of operator-valued functions defined on a subset of $\mathbb{C}$ with values in $\mathcal{L}(\mathcal{B},\mathcal{B}')$ that are holomorphic in a neighborhood of $z_0 \in \C$, except possibly at $z_0$. 

The notions introduced below play roles analogous to eigenvalues, eigenvectors, and algebraic multiplicities in the finite-dimensional theory of matrix pencils.
Consider a function $F \in \mathfrak{U}(z_0)$. We call $z_0$ a \textit{characteristic value} of $F(z)$ if there exists a vector-valued function $f(z)$ with values in $\mathcal{B}$ such that 
    \begin{enumerate}
	\item[(i)] $f(z)$ is holomorphic at $z_0$ and $f(z_0) \neq 0$
	\item[(ii)] $F(z)f(z)$ is holomorphic at $z_0$ and $F(z_0)f(z_0)=0$. 
    \end{enumerate} 
Here, $f(z)$ is called the \textit{root function} of $F(z)$ associated with the characteristic value $z_0$.

Suppose that $z_0$ is a characteristic value of $F(z)$ and $f(z)$ is an associated root function. Then there exists a positive integer $m(f) \in \N$ and a vector-valued function $g$ with values in $\mathcal{B}'$ that is holomorphic at $z_0$ such that  
\[
    F(z)f(z) = (z-z_0)^{m(f)}g(z), \quad \text{with} \ g(z_0) \neq 0.
\] 
The number $m(f)$ is called the \textit{multiplicity} of the root function $f(z)$ and the vector $f(z_0) \in \mathcal{B}$ is called the \textit{eigenvector} of $F(z)$ associated with $z_0$. The closure of all such eigenvectors is the \textit{kernel} of $F(z)$ associated with $z_0$ and denoted by ${\rm Ker}\,F(z_0)$.

For an eigenvector $\phi \in {\rm Ker}\,F(z_0)$, we define the \textit{rank} of $\phi$ by
\[
    {\rm rank}(\phi) = \max \{m(f): f \text{ is a root function with } f(z_0) = \phi \}.
\]
When $n = \dim\, {\rm Ker}\,F(z_0) < \infty$ and ${\rm rank}(\phi) < \infty$ for every $\phi \in {\rm Ker}\,F(z_0)$, the \textit{canonical system of eigenvectors} of $F(z)$ associated with $z_0$ is a collection of eigenvectors $\phi_1,\dots,\phi_n \in {\rm Ker}\,F(z_0)$ chosen recursively as follows. The vector $\phi_1$ is chosen to have maximal rank in ${\rm Ker}\,F(z_0)$. For each $j=2,\dots,n$, the vector $\phi_j$ is chosen to have maximal rank among all eigenvectors belonging to a direct complement of ${\rm span}\{\phi_1,\dots,\phi_{j-1}\}$ in ${\rm Ker}\,F(z_0)$.
The \textit{null multiplicity} of the characteristic value $z_0$ of $F(z)$ is defined by
\[
    N(F(z_0)) = \sum_{j=1}^n {\rm rank}(\phi_j).
\]
We set $N(F(z_0)) =0$ when $z_0$ is not a characteristic value of $F(z)$.

Suppose further that $F^{-1}(z)$ exists and is holomorphic in a deleted neighborhood of $z_0$. If $z_0$ is a characteristic value of $F^{-1}(z)$ and ${\rm dim}\, {\rm Ker}\, F^{-1}(z_0) <\infty$, we define $N(F^{-1}(z_0))$ in the same manner as above. This quantity is called the \textit{polar multiplicity} of $z_0$ for $F(z)$. Then we define the \textit{multiplicity} of $z_0$ for $F(z)$ as
\begin{align*} 
    M(F(z_0)) = N(F(z_0)) - N(F^{-1}(z_0)). 
\end{align*}
It follows that if $z_0$ is a characteristic value of $F(z)$ but not of $F^{-1}(z)$, then $M(F(z_0)) = N(F(z_0))$; and if $z_0$ is a characteristic value of $F^{-1}(z)$ but not of $F(z)$, then $M(F(z_0)) = -N(F^{-1}(z_0))$. When $F(z)$ is holomorphic and invertible at $z_0$, $z_0$ is called a \textit{regular point} of $F(z)$, which results in $M(F(z_0)) = 0$.

Before introducing the remaining concepts, we briefly recall the notion of a Fredholm operator and its index. An operator $A \in L(\mathcal{B},\mathcal{B}')$ is called \textit{Fredholm} provided the subspace $\mathrm{Ker}\, A$ is finite-dimensional and the subspace $\mathrm{Im}\, A$ is closed in $\mathcal{B}'$ and of finite codimension. Let $\mathrm{Fred}(\mathcal{B},\mathcal{B}')$ denote the collection of all Fredholm operators from $\mathcal{B}$ into $\mathcal{B}'$. It is well known that $\mathrm{Fred}(\mathcal{B},\mathcal{B}')$ is an open subset of $\mathcal{L}(\mathcal{B},\mathcal{B}')$.
The \textit{index} of a Fredholm operator $A$ is defined by
\[
    \mathrm{ind}\,A = \mathrm{dim}\, \mathrm{Ker}\, A -  \mathrm{codim}\, \mathrm{Im}\, A. 
\]
    
If $z_0$ is a pole of $F(z)$, then $F(z)$ admits the Laurent series expansion at $z_0$:
\[
    F(z) = \sum_{j= -p}^\infty (z-z_0)^j F_j.
\] 
If $F_{j}$, $-p \leq j \leq -1$, have finite ranks, then $F(z)$ is said to be \textit{finitely meromorphic} at $z_0$. If $F_0$ is a Fredholm operator (of index zero), then $F(z)$ is said to be of \textit{Fredholm type} (of index zero) at $z_0$. If $F(z)$ is finitely meromorphic at $z_0$, is of Fredholm type at $z_0$, and every point in a deleted neighborhood of $z_0$ are regular points of $F(z)$, then $z_0$ is called a \textit{normal point} of $F(z)$. It is known that a normal point of $F(z)$ is also a normal point of $F^{-1}(z)$. 

Consider a simply-connected and bounded domain $D \subset \C$ with rectifiable boundary $\p D$. We say that $F(z)$ is \textit{normal} with respect to $\p D$ if it is finitely meromorphic and of Fredholm type in $D$, continuous on $\p D$, and invertible at all points in $\overline{D}$ except at a finite number of normal points of $F(z)$ in $D$. Suppose that $F(z)$ is normal with respect to $\p D$ and let $z_j$, $j=1,\dots,N$, denote all the characteristic values $F(z)$ and $F^{-1}(z)$ in $D$. Then the \textit{full multiplicity} of $F(z)$ in $D$ is defined by 
\[
    \mathcal{M}(F(z);\p D) = \sum_{j=1}^N M(F(z_j)).
\]
That is, $\mathcal{M}(F(z);\p D)$ equals the number of characteristic values of $F(z)$ in $D$, counted with multiplicities, minus the number of characteristic values of $F^{-1}(z)$ in $D$, counted with multiplicities. 

Under the above assumptions, we conclude this section by stating three fundamental results of the Gohberg--Sigal theory. The generalized argument principle for operator-valued functions takes the following form. 
\begin{thm}[The generalized argument principle {\cite[Theorem 4.1]{Gohberg:1971:OGL}}] \label{thm:GS:argument}
For a scalar function $h$ that is analytic in $D$ and continuous on $\overline{D}$, it holds that
\begin{align*}
    \frac{1}{2\pi i} {\rm tr}\int_{\p D} h(z)F^{-1}(z)\frac{d}{dz}F(z)dz = \sum_{j=1}^N M(F(z_j))h(v_j),
\end{align*} 
where {\rm tr} denotes the trace of an (finite-rank) operator. 
\end{thm}
In particular, when $h \equiv 1$ in the theorem, it holds that
\[
    \frac{1}{2\pi i} {\rm tr}\int_{\p D} F^{-1}(z)\frac{d}{dz}F(z)dz = \mathcal{M}(F(z);\p D).
\]
The contour integrals above are well-defined. Indeed, if $A \in \mathfrak{U}(z_0)$ is finitely meromorphic on $\overline{U}$, where $U$ is a neighborhood of $z_0$ whose closure contains no poles of $A(z)$ except possibly $z_0$, then the operator $\int_{\p U} A(z) dz$ has a finite rank. 

Next, we formulate the generalized version of Rouch\'e theorem. 
\begin{thm}[The generalized Rouch\'{e} theorem {\cite[Theorem 2.2]{Gohberg:1971:OGL}}] \label{thm:GS:Rouche}
Let $G(z)$ be an operator-valued function defined on a neighborhood of $\overline{D} \subset \mathbb{C}$ with values in $\mathcal{L}(\mathcal{B},\mathcal{B}')$. If $G$ is finitely meromorphic in $D$ and continuous on $\p D$ satisfying 
\[
\| F^{-1}(z) G(z)\|_{\mathcal{L}(\mathcal{B},\mathcal{B})} <1, \quad \text{for } z \in \p D, 	\]
then $F(z)+G(z)$ is also normal with respect to $\p D$ and
\[
    \mathcal{M}(F(z)+G(z); \p D) = \mathcal{M}(F(z); \p D).
\]
\end{thm}

Lastly, we introduce the generalized Steinberg theorem. 
\begin{thm}[The generalized Steinberg theorem {\cite[Corollary 3.3]{Gohberg:1971:OGL}}] \label{thm:GS:Steinberg}
Suppose that $F(z)$ is finitely meromorphic and of Fredholm type in $D$. If $F(z)$ is invertible at one point of $D$, then $F(z)$ has a bounded inverse for all $z \in D$, except possibly for some isolated points.  
\end{thm}

\bibliographystyle{plain}
\bibliography{Periodic_ref}

@article{Alouges:2018:FBS,
     author = {Fran\c{c}ois Alouges and Matthieu Aussal},
     title = {FEM and {BEM} simulations with the {Gypsilab} framework},
     journal = {The SMAI Journal of computational mathematics},
     pages = {297--318},
     publisher = {Soci\'et\'e de Math\'ematiques Appliqu\'ees et Industrielles},
     volume = {4},
     year = {2018},
     doi = {10.5802/smai-jcm.36},
     zbl = {1416.65429},
     url = {https://smai-jcm.centre-mersenne.org/articles/10.5802/smai-jcm.36/}
}

@article{Ammari:2004:SRS,
author = {Habib Ammari and Faouzi Triki},
title = {Splitting of resonant and scattering frequencies under shape deformation},
journal = {Journal of Differential Equations},
volume = {202},
number = {2},
pages = {231-255},
year = {2004},
issn = {0022-0396},
doi = {https://doi.org/10.1016/j.jde.2004.02.017},
url = {https://www.sciencedirect.com/science/article/pii/S0022039604001111},
}

@book {Ammari:2009:LPT,
    AUTHOR = {Ammari, Habib and Kang, Hyeonbae and Lee, Hyundae},
     TITLE = {Layer potential techniques in spectral analysis},
    SERIES = {Mathematical Surveys and Monographs},
    VOLUME = {153},
 PUBLISHER = {American Mathematical Society, Providence, RI},
      YEAR = {2009},
     PAGES = {vi+202},
      ISBN = {978-0-8218-4784-8},
       DOI = {10.1090/surv/153},
       URL = {https://doi.org/10.1090/surv/153},
}

@article {Ammari:2015:MTS,
    AUTHOR = {Ammari, Habib and Zhang, Hai},
     TITLE = {A mathematical theory of super-resolution by using a system of
              sub-wavelength {H}elmholtz resonators},
   JOURNAL = {Comm. Math. Phys.},
  FJOURNAL = {Communications in Mathematical Physics},
    VOLUME = {337},
      YEAR = {2015},
    NUMBER = {1},
     PAGES = {379--428},
      ISSN = {0010-3616,1432-0916},
       DOI = {10.1007/s00220-015-2301-4},
       URL = {https://doi.org/10.1007/s00220-015-2301-4}
}

@article{Ammari:2016:MNF,
author = {Ammari, Habib  and Ruiz, Matias  and Wu, Wei  and Yu, Sanghyeon  and Zhang, Hai },
title = {Mathematical and numerical framework for metasurfaces using thin layers of periodically distributed plasmonic nanoparticles},
journal = {Proceedings of the Royal Society A: Mathematical, Physical and Engineering Sciences},
volume = {472},
number = {2193},
pages = {20160445},
year = {2016},
doi = {10.1098/rspa.2016.0445},
URL = {https://royalsocietypublishing.org/doi/abs/10.1098/rspa.2016.0445},
eprint = {https://royalsocietypublishing.org/doi/pdf/10.1098/rspa.2016.0445}
}

@article{Ammari:2017:EMT,
  author  = {Ammari, Habib and Zhang, Hai},
  title   = {Effective Medium Theory for Acoustic Waves
             in Bubbly Fluids near {Minnaert} Resonant Frequency},
  journal = {SIAM Journal on Mathematical Analysis},
  volume  = {49},
  number  = {4},
  pages   = {3252--3276},
  year    = {2017},
  doi     = {10.1137/16M1078574}
}

@article {Ammari:2017:MNF,
    AUTHOR = {Ammari, Habib and Fitzpatrick, Brian and Gontier, David and
              Lee, Hyundae and Zhang, Hai},
     TITLE = {A mathematical and numerical framework for bubble
              meta-screens},
   JOURNAL = {SIAM J. Appl. Math.},
  FJOURNAL = {SIAM Journal on Applied Mathematics},
    VOLUME = {77},
      YEAR = {2017},
    NUMBER = {5},
     PAGES = {1827--1850},
      ISSN = {0036-1399,1095-712X},
       DOI = {10.1137/16M1090235},
       URL = {https://doi.org/10.1137/16M1090235},
}

@article {Ammari:2018:MRA,
    AUTHOR = {Ammari, Habib and Fitzpatrick, Brian and Gontier, David and
              Lee, Hyundae and Zhang, Hai},
     TITLE = {Minnaert resonances for acoustic waves in bubbly media},
   JOURNAL = {Ann. Inst. H. Poincar\'e{} C Anal. Non Lin\'eaire},
  FJOURNAL = {Annales de l'Institut Henri Poincar\'e{} C. Analyse Non
              Lin\'eaire},
    VOLUME = {35},
      YEAR = {2018},
    NUMBER = {7},
     PAGES = {1975--1998},
      ISSN = {0294-1449,1873-1430},
       DOI = {10.1016/j.anihpc.2018.03.007},
       URL = {https://doi.org/10.1016/j.anihpc.2018.03.007},
}

@article {Beceanu:2026:EPL,
    AUTHOR = {Beceanu, Marius and Hong, Jiho and Kwon, Hyun-Kyoung and Lim,
              Mikyoung},
     TITLE = {The eigenvalue problem for the {L}aplacian via conformal
              mapping and the {G}ohberg-{S}igal theory},
   JOURNAL = {Potential Anal.},
  FJOURNAL = {Potential Analysis. An International Journal Devoted to the
              Interactions between Potential Theory, Probability Theory,
              Geometry and Functional Analysis},
    VOLUME = {64},
      YEAR = {2026},
    NUMBER = {3},
     PAGES = {Paper No. 48, 33},
      ISSN = {0926-2601,1572-929X},
       DOI = {10.1007/s11118-026-10283-y},
       URL = {https://doi.org/10.1007/s11118-026-10283-y},
}

@misc{Bocquet:2015:GEI,
  author       = {Stephen Bocquet},
  title        = {Generalised Exponential Integral},
  year         = {2015},
  howpublished = {MATLAB Central File Exchange},
  note         = {Version 1.0.0.0},
  url          = {https://www.mathworks.com/matlabcentral/fileexchange/52694-generalised-exponential-integral}
}

@article{Bruno:2014:RCT,
title = {Rapidly convergent two-dimensional quasi-periodic Green function throughout the spectrum—including Wood anomalies},
journal = {Journal of Computational Physics},
volume = {262},
pages = {262-290},
year = {2014},
issn = {0021-9991},
doi = {https://doi.org/10.1016/j.jcp.2013.12.047},
url = {https://www.sciencedirect.com/science/article/pii/S0021999113008553},
author = {Oscar P. Bruno and Bérangère Delourme}
}

@article{Bruno:2017:RCQ,
author = {Bruno, Oscar P.  and Fernandez-Lado, Agustin G. },
title = {Rapidly convergent quasi-periodic Green functions for scattering by arrays of cylinders—including Wood anomalies},
journal = {Proceedings of the Royal Society A: Mathematical, Physical and Engineering Sciences},
volume = {473},
number = {2199},
pages = {20160802},
year = {2017},
doi = {10.1098/rspa.2016.0802},
URL = {https://royalsocietypublishing.org/doi/abs/10.1098/rspa.2016.0802},
eprint = {https://royalsocietypublishing.org/doi/pdf/10.1098/rspa.2016.0802}
}

@article{Bruno:2017:TDQ,
author = {Bruno, Oscar P.  and Shipman, Stephen P.  and Turc, Catalin  and Stephanos, Venakides },
title = {Three-dimensional quasi-periodic shifted Green function throughout the spectrum, including Wood anomalies},
journal = {Proceedings of the Royal Society A: Mathematical, Physical and Engineering Sciences},
volume = {473},
number = {2207},
pages = {20170242},
year = {2017},
doi = {10.1098/rspa.2017.0242},
URL = {https://royalsocietypublishing.org/doi/abs/10.1098/rspa.2017.0242},
eprint = {https://royalsocietypublishing.org/doi/pdf/10.1098/rspa.2017.0242}
}

@article{Bruno:2020:EQG,
author = {Oscar P. Bruno and Agustin G. Fernandez-Lado},
title = {On the evaluation of quasi-periodic Green functions and wave-scattering at and around Rayleigh-Wood anomalies},
journal = {Journal of Computational Physics},
volume = {410},
pages = {109352},
year = {2020},
issn = {0021-9991},
doi = {https://doi.org/10.1016/j.jcp.2020.109352},
url = {https://www.sciencedirect.com/science/article/pii/S0021999120301261},
}

@article{Calvo:2012:LFS,
    author = {Calvo, David C. and Thangawng, Abel L. and Layman, Christopher N.},
    title = {Low-frequency resonance of an oblate spheroidal cavity in a soft elastic medium},
    journal = {The Journal of the Acoustical Society of America},
    volume = {132},
    number = {1},
    pages = {EL1-EL7},
    year = {2012},
    month = {06},
    issn = {0001-4966},
    doi = {10.1121/1.4721646},
    url = {https://doi.org/10.1121/1.4721646}
}

@article{Capolino:2007:EC3,
author = {F. Capolino and D.R. Wilton and W.A. Johnson},
title = {Efficient computation of the 3D Green’s function for the Helmholtz operator for a linear array of point sources using the Ewald method},
journal = {Journal of Computational Physics},
volume = {223},
number = {1},
pages = {250-261},
year = {2007},
issn = {0021-9991},
doi = {https://doi.org/10.1016/j.jcp.2006.09.013},
url = {https://www.sciencedirect.com/science/article/pii/S0021999106004359},
}

@article{Cho:2019:SAN,
title = {Spectrally-accurate numerical method for acoustic scattering from doubly-periodic 3D multilayered media},
journal = {Journal of Computational Physics},
volume = {393},
pages = {46-58},
year = {2019},
issn = {0021-9991},
doi = {https://doi.org/10.1016/j.jcp.2019.05.010},
url = {https://www.sciencedirect.com/science/article/pii/S0021999119303353},
author = {Min Hyung Cho}
}

@article{Commander:1989:LPW,
    author = {Commander, Kerry W. and Prosperetti, Andrea},
    title = {Linear pressure waves in bubbly liquids: Comparison between theory and experiments},
    journal = {The Journal of the Acoustical Society of America},
    volume = {85},
    number = {2},
    pages = {732-746},
    year = {1989},
    month = {02},
    issn = {0001-4966},
    doi = {10.1121/1.397599},
    url = {https://doi.org/10.1121/1.397599}
}

@Article{Darweesh:2018:RRW,
AUTHOR = {Darweesh, Ahmad A. and Bauman, Stephen J. and Debu, Desalegn T. and Herzog, Joseph B.},
TITLE = {The Role of Rayleigh-Wood Anomalies and Surface Plasmons in Optical Enhancement for Nano-Gratings},
JOURNAL = {Nanomaterials},
VOLUME = {8},
YEAR = {2018},
NUMBER = {10},
ARTICLE-NUMBER = {809},
URL = {https://www.mdpi.com/2079-4991/8/10/809},
ISSN = {2079-4991},
DOI = {10.3390/nano8100809}
}

@article{Errico:2015:UUL,
  author  = {Errico, Claudia and Pierre, Juliette and
             Pezet, Sophie and Desailly, Yann and
             Lenkei, Zsolt and Couture, Olivier and
             Tanter, Mickael},
  title   = {Ultrafast ultrasound localization microscopy
             for deep super-resolution vascular imaging},
  journal = {Nature},
  volume  = {527},
  pages   = {499--502},
  year    = {2015},
  doi     = {10.1038/nature16066}
}

@article{Ewald:1921:BOE,
  title={Die Berechnung optischer und elektrostatischer Gitterpotentiale},
  author={Ewald, Paul P.},
  journal={Annalen der physik},
  volume={369},
  number={3},
  pages={253--287},
  year={1921},
  doi = {https://doi.org/10.1002/andp.19213690304},
  url = {https://onlinelibrary.wiley.com/doi/abs/10.1002/andp.19213690304},
  eprint = {https://onlinelibrary.wiley.com/doi/pdf/10.1002/andp.19213690304},
}

@article {Gohberg:1971:OGL,
    AUTHOR = {Gohberg, I. and Sigal, E.},
     TITLE = {An operator generalization of the logarithmic residue theorem
              and {R}ouch\'e's theorem},
   JOURNAL = {Mat. Sb. (N.S.)},
  FJOURNAL = {Matematicheski\u i\ Sbornik. Novaya Seriya},
    VOLUME = {84(126)},
      YEAR = {1971},
     PAGES = {607--629},
      ISSN = {0368-8666},
}

@article{Hwang:2000:LFR,
  author  = {Hwang, Paul A. and Teague, William J.},
  title   = {Low-Frequency Resonant Scattering of Bubble Clouds},
  journal = {Journal of Atmospheric and Oceanic Technology},
  volume  = {17},
  number  = {6},
  pages   = {847--853},
  year    = {2000},
  doi     = {10.1175/1520-0426(2000)017<0847:LFRSOB>2.0.CO;2}
}

@misc{Johnson:2012:FPC,
  author       = {Johnson, Steven G.},
  title        = {Faddeeva Package: complex error functions},
  year         = {2012},
  howpublished = {MATLAB Central File Exchange},
  note         = {Version 1.5.0},
  url          = {https://www.mathworks.com/matlabcentral/fileexchange/38787-faddeeva-package-complex-error-functions}
}

@article{Lai:2015:FRS,
title = {A fast and robust solver for the scattering from a layered periodic structure containing multi-particle inclusions},
journal = {Journal of Computational Physics},
volume = {298},
pages = {194-208},
year = {2015},
issn = {0021-9991},
doi = {https://doi.org/10.1016/j.jcp.2015.06.005},
url = {https://www.sciencedirect.com/science/article/pii/S0021999115003848},
author = {Jun Lai and Motoki Kobayashi and Alex Barnett}
}

@article{Lanoy:2015:SFB,
  title = {Subwavelength focusing in bubbly media using broadband time reversal},
  author = {Lanoy, Maxime and Pierrat, Romain and Lemoult, Fabrice and Fink, Mathias and Leroy, Valentin and Tourin, Arnaud},
  journal = {Phys. Rev. B},
  volume = {91},
  issue = {22},
  pages = {224202},
  numpages = {6},
  year = {2015},
  month = {Jun},
  publisher = {American Physical Society},
  doi = {10.1103/PhysRevB.91.224202},
  url = {https://link.aps.org/doi/10.1103/PhysRevB.91.224202}
}

@article{Leroy:2009:DCB,
  author  = {Leroy, Valentin and Bretagne, Alice and
             Fink, Mathias and Willaime, Herv{\'e} and
             Tabeling, Patrick and Tourin, Arnaud},
  title   = {Design and Characterization of Bubble Phononic Crystals},
  journal = {Applied Physics Letters},
  volume  = {95},
  number  = {17},
  pages   = {171904},
  year    = {2009},
  doi     = {10.1063/1.3254243}
}

@article{Leroy:2009:TUS,
  author  = {Leroy, Valentin and Strybulevych, Anatoliy and
             Scanlon, Michael G. and Page, John H.},
  title   = {Transmission of Ultrasound through a Single Layer of Bubbles},
  journal = {The European Physical Journal E},
  volume  = {29},
  number  = {1},
  pages   = {123--130},
  year    = {2009},
  doi     = {10.1140/epje/i2009-10457-y}
}

@article{Leroy:2015:SAW,
  title = {Superabsorption of acoustic waves with bubble metascreens},
  author = {Leroy, Valentin and Strybulevych, Anatoliy and Lanoy, Maxime and Lemoult, Fabrice and Tourin, Arnaud and Page, John H.},
  journal = {Phys. Rev. B},
  volume = {91},
  issue = {2},
  pages = {020301(R)},
  numpages = {5},
  year = {2015},
  month = {Jan},
  publisher = {American Physical Society},
  doi = {10.1103/PhysRevB.91.020301},
  url = {https://link.aps.org/doi/10.1103/PhysRevB.91.020301}
}

@article {Linton:2010:LSH,
    AUTHOR = {Linton, Chris M.},
     TITLE = {Lattice sums for the {H}elmholtz equation},
   JOURNAL = {SIAM Rev.},
  FJOURNAL = {SIAM Review},
    VOLUME = {52},
      YEAR = {2010},
    NUMBER = {4},
     PAGES = {630--674},
      ISSN = {0036-1445,1095-7200},
       DOI = {10.1137/09075130X},
       URL = {https://doi.org/10.1137/09075130X},
}

@article{Liu:2016:ENS,
title = {Efficient numerical solution of acoustic scattering from doubly-periodic arrays of axisymmetric objects},
journal = {Journal of Computational Physics},
volume = {324},
pages = {226-245},
year = {2016},
issn = {0021-9991},
doi = {https://doi.org/10.1016/j.jcp.2016.08.011},
url = {https://www.sciencedirect.com/science/article/pii/S0021999116303540},
author = {Yuxiang Liu and Alex H. Barnett}
}

@article{Maradudin:2016:RWA,
doi = {10.1088/2040-8978/18/2/024004},
url = {https://doi.org/10.1088/2040-8978/18/2/024004},
year = {2016},
month = {jan},
publisher = {IOP Publishing},
volume = {18},
number = {2},
pages = {024004},
author = {Maradudin, A A and Simonsen, I and Polanco, J and Fitzgerald, R M},
title = {Rayleigh and Wood anomalies in the diffraction of light from a perfectly conducting reflection grating},
journal = {Journal of Optics}
}

@Inbook{Maystre:2012:TWA,
author= {Maystre, Daniel},
editor={Enoch, Stefan and Bonod, Nicolas},
title={Theory of Wood's Anomalies},
booktitle={Plasmonics: From Basics to Advanced Topics},
year={2012},
publisher={Springer Berlin Heidelberg},
address={Berlin, Heidelberg},
pages={39--83},
isbn={978-3-642-28079-5},
doi={10.1007/978-3-642-28079-5_2},
url={https://doi.org/10.1007/978-3-642-28079-5_2}
}

@article{Minnaert:1933:MAB,
author = {M. Minnaert},
title = {XVI. On musical air-bubbles and the sounds of running water },
journal = {The London, Edinburgh, and Dublin Philosophical Magazine and Journal of Science},
volume = {16},
number = {104},
pages = {235--248},
year = {1933},
publisher = {Taylor \& Francis},
doi = {10.1080/14786443309462277},
URL = {https://doi.org/10.1080/14786443309462277}
}

@article{Nicholls:2020:SPI,
  author  = {Nicholls, David P. and
             P{\'e}rez-Arancibia, Carlos and
             Turc, Catalin},
  title   = {Sweeping Preconditioners for the Iterative
             Solution of Quasiperiodic {Helmholtz}
             Transmission Problems in Layered Media},
  journal = {Journal of Scientific Computing},
  volume  = {82},
  number  = {44},
  year    = {2020},
  doi     = {10.1007/s10915-020-01133-z}
}

@article{Perez-Arancibia:2019:DDQ,
  author = {Pérez-Arancibia, Carlos and Shipman, Stephen P. and Turc, Catalin and Venakides, Stephanos},
  title = {Domain Decomposition for Quasi-Periodic Scattering by Layered Media via Robust Boundary-Integral Equations at All Frequencies},
  journal = {Communications in Computational Physics},
  year = {2019},
  volume = {26},
  number = {1},
  pages = {265-310},
  doi = {10.4208/cicp.OA-2018-0021},
  url = {https://doi.org/10.4208/cicp.OA-2018-0021}
}

@article {Petit:1990:ETG,
    AUTHOR = {Petit, R. and Cadilhac, M.},
     TITLE = {Electromagnetic theory of gratings: some advances and some
              comments on the use of the operator formalism},
   JOURNAL = {J. Opt. Soc. Amer. A},
  FJOURNAL = {Journal of the Optical Society of America A},
    VOLUME = {7},
      YEAR = {1990},
    NUMBER = {9},
     PAGES = {1666--1674},
      ISSN = {0740-3232},
   MRCLASS = {78A45 (78-02)},
  MRNUMBER = {1068227},
       DOI = {10.1364/JOSAA.7.001666},
       URL = {https://doi.org/10.1364/JOSAA.7.001666},
}

@article{Savoia:2013:SSR,
  author  = {Savoia, Silvio and Ricciardi, Armando and
             Crescitelli, Alessio and Granata, Carmine and
             Esposito, Emanuela and Galdi, Vincenzo and
             Cusano, Andrea},
  title   = {Surface sensitivity of {Rayleigh} anomalies in
             metallic nanogratings},
  journal = {Optics Express},
  volume  = {21},
  number  = {20},
  pages   = {23531--23542},
  year    = {2013},
  doi     = {10.1364/OE.21.023531}
}

@article{Tavakoli:2019:RWA,
title = {Rayleigh-Wood anomaly approximation with FDTD simulation of plasmonic gold nanohole array for determination of optimum extraordinary optical transmission characteristics},
journal = {Superlattices and Microstructures},
volume = {130},
pages = {454-471},
year = {2019},
issn = {0749-6036},
doi = {https://doi.org/10.1016/j.spmi.2019.04.035},
url = {https://www.sciencedirect.com/science/article/pii/S0749603618324820},
author = {Mehdi Tavakoli and Yousef Seyed Jalili and Seyed Mohammad Elahi}
}

@article{Valerio:2007:CAA,
  author={Valerio, Guido and Baccarelli, Paolo and Burghignoli, Paolo and Galli, Alessandro},
  journal={IEEE Transactions on Antennas and Propagation}, 
  title={Comparative Analysis of Acceleration Techniques for 2-D and 3-D Green's Functions in Periodic Structures Along One and Two Directions}, 
  year={2007},
  volume={55},
  number={6},
  pages={1630-1643},
  doi={10.1109/TAP.2007.897340}}

\end{document}